\documentclass[11pt]{amsart}
\usepackage[margin=1in]{geometry}
\usepackage[utf8]{inputenc}
\usepackage{amsthm, amssymb}
\usepackage[colorlinks=true, pdfstartview=FitV, linkcolor=blue, citecolor=blue, urlcolor=blue]{hyperref}
\usepackage[colorinlistoftodos]{todonotes}
\usepackage{algorithm}
\usepackage{algpseudocode}
\usepackage{pgf}
\usepackage{tikz}
\usepackage{tikz-cd}
\usetikzlibrary{positioning,shapes,shadows,arrows}
\usepackage{bbm,bm}
\usepackage{enumitem}
\usepackage{mathabx}
\usepackage{comment}
\usepackage{ytableau}
\usepackage{thmtools}
\usepackage{thm-restate}

\newtheorem{prop}{Proposition}[section]
\newtheorem{thm}[prop]{Theorem}

\newtheorem{lemma}[prop]{Lemma}
\newtheorem{cor}[prop]{Corollary}
\newtheorem{conj}[prop]{Conjecture}
\newtheorem*{thm_main}{Theorem~\ref{T: FS on BPD}}
\newtheorem*{thm_main_chain}{Theorem~\ref{T: BPD to chain}}
\newtheorem*{thm_main_bijection}{Theorem~\ref{T: BPD to PD bijection}}
\newtheorem*{cor_Bruhat}{Corollary~\ref{C: Bruhat formula}}
\theoremstyle{remark}
\newtheorem{exa}[prop]{Example}
\newtheorem{rem}[prop]{Remark}
\newtheorem{defn}[prop]{Definition}
\newtheorem{prob}[prop]{Problem}

\newcommand{\E}{\mathcal{E}}
\newcommand{\N}{\mathcal{N}}

\newcommand{\BPD}{\mathsf{BPD}}
\newcommand{\HPD}{\mathsf{HPD}}

\newcommand{\fS}{\mathfrak{S}}

\newcommand{\chain}{\mathsf{chain}}
\newcommand{\PD}{\mathsf{PD}}

\newcommand{\word}{\mathsf{word}}
\newcommand{\flip}{\mathsf{flip}}
\newcommand{\wt}{\mathsf{wt}}
\newcommand{\WT}{\mathsf{WT}}
\newcommand{\Z}{\mathbb{Z}}
\newcommand{\Q}{\mathbb{Q}}

\newcommand{\FT}{\mathsf{FT}}
\newcommand{\growth}{\mathsf{growth}}

\newcommand{\btile}{
 \begin{tikzpicture}[x=1em,y=1em,thick,color = blue]
\draw[step=1,gray,thin] (0,0) grid (1,1);
\draw[color=black, thick, sharp corners] (0,0) rectangle (1,1);
\end{tikzpicture}}

\newcommand{\htile}{
\begin{tikzpicture}[x=1em,y=1em,thick,color = blue]
\draw[step=1,gray,thin] (0,0) grid (1,1);
\draw[color=black, thick, sharp corners] (0,0) rectangle (1,1);
\draw(1.0,0.5)--(0.0,0.5);
\end{tikzpicture}
}

\newcommand{\vtile}{
\begin{tikzpicture}[x=1em,y=1em,thick,color = blue]
\draw[step=1,gray,thin] (0,0) grid (1,1);
\draw[color=black, thick, sharp corners] (0,0) rectangle (1,1);
\draw(0.5,1.0)--(0.5,0.0);
\end{tikzpicture}}

\newcommand{\ptile}{
\begin{tikzpicture}[x=1em,y=1em,thick,color = blue]
\draw[step=1,gray,thin] (0,0) grid (1,1);
\draw[color=black, thick, sharp corners] (0,0) rectangle (1,1);
\draw(0.5,1.0)--(0.5,0.0);
\draw(1.0,0.5)--(0.0,0.5);
\end{tikzpicture}}

\newcommand{\rtile}{
\begin{tikzpicture}[x=1em,y=1em,thick,rounded corners,color = blue]
\draw[step=1,gray,thin] (0,0) grid (1,1);
\draw[color=black, thick, sharp corners] (0,0) rectangle (1,1);
\draw(1.0,0.5)--(0.5,0.5)--(0.5,0.0);
\end{tikzpicture}}

\newcommand{\jtile}{
\begin{tikzpicture}[x=1em,y=1em,thick,rounded corners,color = blue]
\draw[step=1,gray,thin] (0,0) grid (1,1);
\draw[color=black, thick, sharp corners] (0,0) rectangle (1,1);
\draw(0.5,1.0)--(0.5,0.5)--(0.0,0.5);
\end{tikzpicture}}

\newcommand{\Ltile}{
\begin{tikzpicture}[x=1em,y=1em,thick,rounded corners,color = blue]
\draw[step=1,gray,thin] (0,0) grid (1,1);
\draw[color=black, thick, sharp corners] (0,0) rectangle (1,1);
\draw(1.0,0.5)--(0.5,0.5)--(0.5,1.0);
\end{tikzpicture}}

\newcommand{\qtile}{
\begin{tikzpicture}[x=1em,y=1em,thick,rounded corners,color = blue]
\draw[step=1,gray,thin] (0,0) grid (1,1);
\draw[color=black, thick, sharp corners] (0,0) rectangle (1,1);
\draw(0.5,0.0)--(0.5,0.5)--(0.0,0.5);
\end{tikzpicture}}

\newcommand{\bumptile}{
\begin{tikzpicture}[x=1em,y=1em,thick,rounded corners,color = blue]
\draw[step=1,gray,thin] (0,0) grid (1,1);
\draw[color=black, thick, sharp corners] (0,0) rectangle (1,1);
\draw(1.0,0.5)--(0.5,0.5)--(0.5,0.0);
\draw(0.5,1.0)--(0.5,0.5)--(0.0,0.5);
\end{tikzpicture}}

\definecolor{darkblue}{rgb}{0.0,0,0.7} % darkblue color
\definecolor{darkred}{rgb}{0.7,0,0} % darkred color
\definecolor{darkgreen}{rgb}{0, .6, 0} % darkgreen color
\newcommand{\definition}[1]{{\color{darkred}\emph{#1}}} % emphasis of a definition

\title{Embedding bumpless pipedreams as Bruhat chains}
\author[T.~Yu]{Tianyi Yu}
\address[T. Yu]{Department of Mathematics, UC San Diego, La Jolla, CA 92093, U.S.A.}
\email{tiy059@ucsd.edu}

\begin{document}

\begin{abstract}
Schubert polynomials are distinguished representatives of Schubert cycles in the cohomology of the flag variety. 
In the spirit of Bergeron and Sottile, 
we use the Bruhat order to give $(n-1)!$ different combinatorial formulas
for the Schubert polynomial of a permutation in $S_n$.
By work of Lenart and Sottile, 
one extreme of the formulas recover the classical Pipedream (PD) formula.
We prove the other extreme corresponds to Bumpless pipedreams (BPDs).
We give two applications of this perspective to view BPDs: 
Using the Fomin-Kirrilov algebra, 
we solve the problem of finding a BPD analogue
of Fomin and Stanley's algebraic construction on PDs;
We also establish a bijection between PDs and BPDs using Lenart's growth diagram,
which conjectually agrees with the 
existing bijection of Gao and Huang. 
\end{abstract}

\maketitle
\section{Introduction}
Fix $n \in \mathbb{Z}_{\geq 0}$. 
For a permutation $w \in S_n$,
Lascoux and Sch\"utzenberger~\cite{LS:Schubert}
recursively define the \definition{Schubert polynomial} 
$\fS_w$.
The base case is 
$\fS_{w_0} := x_{1}^{n-1} x_2^{n-2} \cdots x_{n-1}$
where $w_0$ is the permutation with 
one-line notation $[n, n-1, \cdots, 1]$.
To compute $\fS_w$ for other $w \in S_n$,
we need the \definition{divided difference operator} 
$\partial_i(f) := \frac{f - f(\cdots, x_{i+1}, x_i, \cdots)}{x_i - x_{i+1}}$.
For $1 \leq i < j \leq n$, we use $t_{i,j} \in S_n$
to denote the transposition that swaps $i$ and $j$.
Then for any $w \in S_n$ and $i \in [n-1]$:
$$
\partial_i(\fS_w) = \begin{cases}
\fS_{w t_{i, i+1}} & \text{if $w(i) > w(i+1)$,} \\
0  & \text{if $w(i) < w(i+1)$.}
\end{cases}
$$

The Schubert polynomials represent Schubert cycles in flag varieties and have been extensively investigated.
Schubert polynomials have two distinct combinatorial formulas involving ``pipes'': pipedreams (PD)~\cite{BB, BJS} and bumpless pipedreams (BPD)~\cite{LLS}.
Both are fillings of grids with certain tiles. 
When we refer to cells of a grid, 
we use the matrix coordinates:
row $1$ is the topmost row 
and column $1$ is the leftmost column. 
A \definition{pipedream} is a filling of a staircase
grid:
The grid has a cell in row $i$ column $j$ for each $i + j \leq n + 1$.
The rightmost cell in each row is $\jtile$.
The rest of the cells can be $\ptile$ (crossing)
or $\bumptile$ (bump), but two pipes cannot cross more than once.
We imagine the pipes enter the grid from the 
left edge and exit from the top edge. 
A \definition{bumpless pipedream (BPD)}
is a consistent filling of an $n \times n$ grid with six types of cells:
$
\vtile, \htile, \ptile, \jtile, \rtile$
and $\btile$.
The pipes enter the grid from each cell
on the right edge and exits on the bottom edge. 
In addition, two pipes cannot cross more than once.
We attach a permutation to each PD (resp. BPD).
The permutation $w$ associated to each PD (resp. BPD) 
can be read off as follows:
For each $i \in [n]$,
let $w(i)$ be the column index of the tile
where the pipe from row $i$ exits.

\begin{exa}
\label{E: PD}
When $n = 5$, we present a PD
and a BPD associated with $[2,5,1,4,3]$:

$$
\begin{tikzpicture}[x=1.5em,y=1.5em,thick,rounded corners, color = blue]
\draw[step=1,black, thin] (0,0) grid (1,1);
\draw[step=1,black,thin] (0,1) grid (2,2);
\draw[step=1,black,thin] (0,2) grid (3,3);
\draw[step=1,black,thin] (0,3) grid (4,4);
\draw[step=1,black,thin] (0,4) grid (5,5);
\draw(0, 0.5)--(0.5,0.5)--(0.5,2.5)--(1.5,2.5)--(1.5,4.5)--(2.5,4.5)--(2.5,5);
\draw(0, 1.5)--(1.5,1.5)--(1.5,2.5)--(2.5,2.5)--(2.5,3.5)--(3.5,3.5)--(3.5,5);
\draw(0, 2.5)--(0.5,2.5)--(0.5,5);
\draw(0, 3.5)--(2.5,3.5)--(2.5,4.5)--(4.5,4.5)--(4.5,5);
\draw(0, 4.5)--(1.5,4.5)--(1.5,5);
\end{tikzpicture}
\quad\quad\quad
\begin{tikzpicture}[x=1.5em,y=1.5em,thick,rounded corners, color = blue]
\draw[step=1,gray,thin] (0,0) grid (5,5);
\draw[color=black, thick, sharp corners] (0,0) rectangle (5,5);
\draw(0.5, 0)--(0.5,2.5)--(5,2.5);
\draw(1.5, 0)--(1.5,3.5)--(2.5,3.5)--(2.5,4.5)--(5,4.5);
\draw(2.5, 0)--(2.5,0.5)--(5,0.5);
\draw(3.5, 0)--(3.5,1.5)--(5,1.5);
\draw(4.5, 0)--(4.5,3.5)--(5,3.5);
\end{tikzpicture}
$$

\end{exa}

Let $\PD(w)$ (resp. $\BPD(w)$) be the set of all
PDs (resp. BPDs) associated 
with $w \in S_n$.
For $P \in \PD(w)$ (resp. $P \in \BPD(w)$),
the weight of $P$, denoted as $\wt(P)$, 
is a sequence of $n-1$ integers where the $i^{th}$ entry
is the number of $\ptile$ (resp. $\btile$) on row $i$.
For instance, the PD and the BPD in Example~\ref{E: PD}
both have weight $(2,2,0,1)$.
If $\alpha = (\alpha_1, \cdots, \alpha_{n-1})$ 
is a sequence of $n-1$ non-negative integers,
we use $x^{\alpha}$ to denote the monomial
$x_1^{\alpha_1} \cdots x_{n-1}^{\alpha_{n-1}}$.
For $w \in S_n$, 
one may compute $\fS_w$
as $\sum_{P \in \PD(w)} x^{\wt(P)}$ by~\cite{BB, BJS}
and as $\sum_{D \in \BPD(w)} x^{\wt(D)}$ by~\cite{LLS}.

Another combinatorial formula of Schubert polynomials
is the Bruhat chain formula of Bergeron and Sottile~\cite{BS}. 
For $w \in S_n$,
let $\ell(w) := |\{(i,j): i < j, w(i) > w(j)|$.
Let $\leq$ be the \definition{Bruhat order}
on $S_n$, 
where the cover relation is 
$u \lessdot w \textrm{ if $w = ut_{i,j}$ and $\ell(w) = \ell(u) + 1$}$.
We say $C = (w_1, w_2, \cdots, w_d)$ 
is a \definition{Bruhat chain}
from $w_1$ to $w_d$ if $w_1 \leq w_2 \leq \cdots \leq w_d$.
The \definition{weight} of $C$, denoted as $\wt(C)$, 
is a sequence of length $d-1$ 
where the $i^{th}$ entry
is $\ell(w_{i+1}) - \ell(w_i)$.
Take $k \in [n-1]$.
We let $\leq_k$ be the
\definition{$k$-Bruhat order}
whose cover relation is given by
$u \lessdot_k w$ if $u \lessdot w$ and 
$w = ut_{i,j}$ for some 
$i \leq k < j$.
We say $(w_1, \cdots, w_d)$ is an \definition{increasing 
$k$-chain} if $w_1 \lessdot_k \cdots \lessdot_k w_d$
and the smaller value being swapped increases from left to right.

We introduce a notion to describe
Bergeron and Sottile's formula in a broader setting. 
We say a Bruhat chain $(w_1, w_2, \cdots, w_{d})$
is \definition{compatible} with a sequence $(k_1, \cdots, k_{d-1})$ if 
there exists an increasing $k_i$-chain from $w_i$
to $w_{i+1}$ for each $i \in [d-1]$.
Then by~\cite{BS},
\begin{equation}
\label{EQ: Bruhat chains}
\fS_w = \sum_{C} x^{(n-1, \cdots, 1) - \wt(C)},
\end{equation}
where the sum is over all chains
from $w$ to $w_0$ compatible with $(1, 2, \cdots, n-1)$.
This formula is closely related to the PD formula:
Lenart and Sottile~\cite{LeSo} biject PDs and chains in~\eqref{EQ: Bruhat chains}. 
(see the map $\chain_\PD$ described in \S~\ref{S: Background}).

To unravel the connection between BPDs and Bruhat chains, 
we first use Sottile's Pieri formula for Schubert polynomials~\cite{Sot}
to extend~\eqref{EQ: Bruhat chains}.
\begin{cor_Bruhat}
Take $w \in S_n$ and $\gamma \in S_{n-1}$.
Then 
$\fS_w = \sum_{C} x^{\wt_\gamma(C)}$,
summing over all chains $C$ from $w$ 
to $w_0$ compatible with $\gamma$.
The summand $x^{\wt_\gamma(C)}$
is a monomial where the power of $x_i$
equals $n - i$ minus the length
of the increasing $i$-chain in $C$.
\end{cor_Bruhat}

In other words, the Bruhat chain formulas of 
Schubert polynomials come in $(n-1)!$ different flavors,
corresponding to $(n-1)!$ different choices of $\gamma \in S_{n-1}$.
Under the map $\chain_\PD$ of Lenart and Sottile,
PDs correspond to an extreme flavor: $\gamma = [1, \cdots, n-1]$.
We will define the map
$\chain_\BPD$ that sends BPDs to chains on the other
extreme: $\gamma = [n-1, \cdots, 1]$.
Our description of $\chain_\BPD$ 
resembles the description of $\chain_\PD$.
Then our first main result is the following: 
\begin{thm_main_chain}
Take $w \in S_n$.
The map $\chain_\BPD$
is a bijection from $\BPD(w)$
to chains from $w$ to $w_0$
that is compatible with $(n-1, \cdots, 1)$.
Moreover, if $\chain_\BPD(D) = C$,
then $\wt(C)$ is the reverse 
of $(n-1, \cdots, 1) - \wt(D)$.
\end{thm_main_chain}

We then provide two applications
of this embedding.
There is a recent surge of research
connecting BPDs with PDs and finding 
BPD analogue of classical PD apparatus
\cite{GH, KU, KW, W}.
Since the introduction of BPDs, 
finding a BPD analogue of the Fomin-Stanley 
construction~\cite{FS} has been an open problem.
Let $u_1, \cdots, u_{n-1}$ be the standard generators
of the nil-Coexter algebra $\N_n$.
Fomin and Stanley defined the following elements 
in $\mathbb{Q}[x_1, \cdots, x_{n-1}] \otimes \N_n$:
$$
A_i(x_i) := (1 + x_i u_{n-1}) (1 + x_i u_{n-2}) \cdots (1 + x_i u_{i}) \textrm{ and }
\fS^\PD := A_1(x_1) \cdots A_{n-1}(x_{n-1}).
$$
Combinatorially, after expanding $\fS^\PD$,
each term $x^\alpha u_{i_1}\cdots u_{i_k}$ 
naturally corresponds to a $P \in \PD(w)$
with $\alpha = \wt(P)$ and $i_1\cdots i_k$ is a reduced word of $w$.
Algebraically, Fomin and Stanley
proved $\fS^\PD = \sum_{w \in S_n} \fS_w u_{i_1}\cdots u_{i_l}$
where $i_1 \cdots i_l$ 
is any reduced word of $w$.
Consequently, they obtain an operator theoretic proof of the PD formula.

We provide a BPD analogue of the Fomin-Stanley construction. 
We consider the Fomin-Kirillov 
algebra $\E_n$~\cite{FK-algebra}.
It is generated by $d_{i,j}$ 
for $1 \leq i < j \leq n$
and has a right action on $\mathbb{Q}[S_n]$:
For $w \in S_n$, $w \odot d_{i,j}:= wt_{i,j}$
if $wt_{i,j} \lessdot w$ and $0$ otherwise.
We then define the following elements 
in $\mathbb{Q}[x_1, \cdots, x_{n-1}] \otimes \E_n$:
$$
R_i(x_i) := 
(x_i + d_{1, i+1} + \cdots + d_{i,i+1})
(x_i + d_{1, i+2} + \cdots + d_{i,i+2})
\cdots
(x_i + d_{1, n} + \cdots + d_{i,n}), \textrm{ and}
$$
$$
\fS^\BPD := w_0 \odot (R_1(x_1) R_2(x_2) \cdots R_{n-1}(x_{n-1})).
$$
Combinatorially, after expanding $\fS^\BPD$,
we show each term $x^\alpha w$
naturally corresponds to a $D \in \BPD(w)$
with $\alpha = \wt(D)$
using Theorem~\ref{T: BPD to chain}.
Algebraically, we establish Theorem~\ref{T: FS on BPD}, obtaining an operator
theoretic proof of the BPD formula.
\begin{thm_main}
We have $\fS^\BPD = \sum_{w \in S_n} \fS_w w$.
\end{thm_main}

Another application of Theorem~\ref{T: BPD to chain}
gives a new perspective on bijecting PDs and BPDs. 
When BPDs were introduced, 
finding a bijection between PDs and BPDs 
became a problem that received considerable attention.
The first solution was given by Gao and Huang~\cite{GH}.
Their bijection relies on repetitively applying
their insertion-like algorithm
on BPDs called ``push''. 
They also showed their bijection
is ``canonical'' in the sense 
of preserving the Monk's rule.
Recently, Knutson and Udell~\cite{KU}
gave another bijection. 
They introduced Hybrid pipedreams (HPDs),  
which are fillings of an $n \times n$
grid and each row may chosen to be
a ``PD row'' or a ``BPD row''.
These objects recover PDs (resp. BPDs)
when all rows are PDs (resp. BPDs) rows.
They defined two local moves 
that can gradually change PD rows
into BPD rows.
Using these moves, 
they provide another bijection from PDs to BPDs 
and verified it agrees with the Gao-Huang bijection.

We provide yet another bijection between PDs and BPDs. 
Thanks to the maps $\chain_\PD$ and $\chain_\BPD$,
this problem reduces to finding a bijection 
between chains from $w$ to $w_0$ 
that are compatible with $(1, \cdots, n-1)$ 
and those compatible with $(n-1, \cdots, 1)$.
We define a map $\flip$ in 
Definition~\ref{D: BPD PD Bijection}
to achieve this objective.
\begin{thm_main_bijection}
The map $\flip$
is a weight-reversing bijection
between chains from $w$ to $w_0$ compatible with $(n-1, \cdots, 1)$
and chains from $w$ to $w_0$ compatible with $(1, \cdots, n-1)$.
Consequently, together with Theorem~\ref{T: BPD to chain} 
and~\cite{LeSo},
it can be viewed as a bijection from $\BPD(w)$ to $\PD(w)$.
\end{thm_main_bijection}
In contrast to the two existing bijections
involving novel algorithms on BPDs or HPDs, 
our bijection $\flip$ uses Lenart's growth diagram~\cite{L_growth},
a classical operation on Bruhat chains.
We conjecture our bijection
agrees with the existing bijections
and have verified this conjecture up to $S_7$.

To summarize, 
PDs and BPDs correspond to two extremes of 
a family of $(n-1)!$ formulas for Schubert polynomials 
in Corollary~\ref{C: Bruhat formula}.
How about the other formulas in this family?
Do they also correspond to combinatorial 
objects involving pipes? 
We can partially answer this question
using the HPDs of Knutson and Udell~\cite{KU}.
In Theorem~\ref{T: HPD to chain},
we show that they correspond 
to $2^{n-2}$ many Bruhat chain formulas from 
Corollary~\ref{C: Bruhat formula}.
We believe it would be an interesting problem
to find a pipe-like formula for arbitrary $\gamma$. 
Besides, 
we also hope various studies on PDs and BPDs
can be extended to the $(n-1)!$ formulas 
(see \S\ref{S: Future} for a list of such questions.)

\textbf{Organization}:
In~\S\ref{S: Background},
we cover some necessary background
and prove Corollary~\ref{C: Bruhat formula}.
In~\S\ref{S: Chain}, 
we describe the map 
$\chain_\BPD$ and prove Theorem~\ref{T: BPD to chain}.
Our first application is
the BPD analogue of Fomin-Stanley construction
in~\S\ref{S: FS}.
Our second application is
a bijection between PDs and BPDs 
using Lenart's growth diagram in~\S\ref{S: Growth}.
Then in~\S\ref{S: HPD},
we give an extension of $\chain_\PD$ 
and $\chain_\BPD$ that biject
HPDs with certain Bruhat chains. 
Finally, in~\S\ref{S: Future}, 
we discuss some future directions.

\section{Background}
\label{S: Background}

\subsection{Lenart and Sottile's bijection}
Lenart and Sottile~\cite{LeSo} described a bijection from 
$\PD(w)$ to chains from $w$ to $w_0$ 
that is compatible with $(1, 2, \cdots, n-1)$.
We denote this map as $\chain_\PD$
and describe this map in an equivalent way.
Our description will parallel
the analogous map on BPDs introduced in \S\ref{S: Chain}. 
Take $P \in \PD(w)$.
For $i \in [n]$,
let $P_i$ be the pipedream
obtained from $P$ by keeping the 
last $i$ rows of $P$.
Let $u_i' \in S_{i}$ be the permutation associated with
$P_i$.
Then we append $n, n-1 \cdots, i+1$ in front
of $u_i'$ to obtain $u_i \in S_n$.
Define $\chain_\PD(P) = (u_n, \cdots, u_{1})$.
Lenart and Sottile showed that
there is an increasing $(n + 1 - k)$-chain
from $u_{k}$ to $u_{k+1}$.
More explicitly, 
let $c_1< \cdots< c_m$ be column indices of $\bumptile$
in row $(n-k+1)$ of $P$,
then the increasing $(n + 1 - k)$-chain from $u_{k}$ to $u_{k-1}$
is given by 
$$
u_{k} \xrightarrow{t_{n-k+1, u_{k-1}^{-1}(c_1)}} \cdots 
\xrightarrow{t_{n-k+1, u_{k-1}^{-1}(c_m)}} u_{k-1}.
$$
In other words, 
$u_{k}$ is obtained from $u_{k-1}$ by swapping the values 
$c_m, \cdots, c_1$ with the $(n-k+1)\textsuperscript{th}$ entry. 

\begin{exa}
\label{E: Lenart Sottile}  
Let $P$ be the pipedream in Example~\ref{E: PD}.
Then
$$
\chain_\PD(P) = 
([2,5,1,4,3], [5,3,1,4,2], [5,4,1,3,2],[5,4,3,2,1],[5,4,3,2,1])
$$
which is compatible with $(1,2,3,4)$.
If we want to obtain the increasing $1$-chain
from $[2,5,1,4,3]$ to $[5,3,1,4,2]$,
we look at column indices of $\bumptile$
in row $1$: They are $2$ and $3$.
Therefore, starting from $[5,3,1,4,2]$,
we first swap the value $3$ with the 
first entry.
Then we swap the value $2$ with the first entry.
The resulting increasing $1$-chain is:
$$
[2,5,1,4,3]
\xrightarrow{\quad\quad t_{1,5} \quad\quad\:} [3,5,1,4,2]
\xrightarrow{\quad\quad t_{1,2} \quad\quad\:} [5,3,1,4,2].
$$
\end{exa}

\begin{thm}\cite{LeSo}
\label{T: PD to chain}
Take $w \in S_n$.
The map $\chain_{\PD}$
is a bijection from $\PD(w)$
to chains from $w$ to $w_0$
that is compatible with $(1, \cdots, n-1)$.
Moreover, if $\chain_\BPD(D) = C$,
then $\wt(C) = (n-1, \cdots, 1) - \wt(D)$.
\end{thm}

We make some observations that will be useful later.
\begin{rem}
\label{R: wt of inc}
Let $(w_1, \cdots, w_d)$ be an increasing $k$-chain
with $w_i = w_{i+1} t_{a_i, b_i}$.
Clearly, $b_1, \cdots, b_{d-1}$ must be distinct.
Observe that for $j > k$,
$w_1(j) \neq w_d(j)$ if and only if
$j$ is among $b_1, \cdots, b_{d-1}$.
Thus, the length of this $k$-chain, 
which is $d-1$,
agrees with $n-k$ subtract the number 
of $j > k$ such that $w_1(j) = w_d(j)$.
\end{rem}

\begin{rem}
\label{R: Double PD}
Suppose $\chain_\PD(P) = (u_n, \cdots, u_1)$.
In $u_{k-1}$, 
the last $k-1$ numbers is a permutation of $[k-1]$.
For a number $j$ among these numbers,
the value $j$ is swapped in the $(n-k+1)$-increasing
chain from $u_k$ to $u_{k-1}$
if and only if $(n-k+1, j)$ is a $\bumptile$.
Following Remark~\ref{R: wt of inc}, 
a number $j$ is among the last 
$k-1$ numbers of $u_{k-1}$ and remains 
unswapped in the $(n-k+1)$-increasing
chain from $u_k$ to $u_{k-1}$ 
if and only if $(n-k+1, j)$ is $\ptile$.
\end{rem}

\subsection{Various Bruhat chain formulas for Schubert polynomials}
We will obtain $(n-1)!$ Bruhat chain formulas
for $\fS_w$ in Corollary~\ref{C: Bruhat formula}.
The strategy is to extend~\eqref{EQ: Bruhat chains}
of Bergeron and Sottile using Sottile's Pieri rule. 
Recall that the $k$-Bruhat order can be used to describe
the Monk's rule~\cite{Monk}:
\begin{equation}
\label{eq: Monk}
\fS_w (x_1 + \cdots + x_k) = \sum_{w \lessdot_k u} \fS_u
\end{equation}
for any $w \in S_n$ and $k \in [n-1]$
such that $w(j) = n$ for some
$j > k$.
Sottile generalized the Monk's rule 
by considering multiplying $\fS_w$
with 
$$
h_d(x_1, \cdots, x_k) := \sum_{1 \leq i_1 \leq \cdots \leq i_d \leq k}x_{i_1} \cdots x_{i_d},
$$
where $k \in [n-1]$ and $d \in \mathbb{Z}_{> 0}$.

\begin{thm}\cite{Sot}
\label{T: Pieri}
Take $u \in S_n$ and $d \in \Z_{\geq 0}$.
For any $k \in [n-1]$ such that
$n - d + 1 \cdots, n$ are among $w(k+1), \cdots, w(n)$ in increasing order,
then
$$
\fS_u \times h_d(x_1, \cdots, x_k) = \sum_{w} \fS_w.
$$
The sum is over all $w$ such that 
$\ell(w) - \ell(u) = d$ and
there is an increasing $k$-chain 
from $u$ to $w$.
\end{thm}

A Bruhat chain $(w_1, \cdots, w_d)$ 
is \definition{saturated} if $w_i \lessdot w_{i+1}$
for each $i$.
We may represent a saturated chain as
$$
w_1 \xrightarrow{t_{a_1, b_1}} w_2 \xrightarrow{t_{a_2, b_2}}
\cdots \xrightarrow{t_{a_{d-1}, b_{d-1}}} w_d,
$$
where $a_i < b_i$ and $w_{i+1} = w_i t_{a_i, b_i}$.
We make one observation regarding increasing $k$-chains. 
\begin{rem}
\label{R: Pieri and increasing chains}
Let $u, w \in S_n$ and $k \in [n-1]$. 
There is at most one increasing 
$k$-chain from $u$ to $w$.
Moreover, if $u \xrightarrow{t_{a_1, b_1}}
\cdots \xrightarrow{t_{a_{d-1}, b_{d-1}}} w$
is increasing, the numbers $b_1,\cdots, b_{d-1}$
must be distinct. 
\end{rem}

\begin{exa}
\label{E: increasing chain}
Suppose $n= 6$.
The unique increasing $2$-chain from 
$[3,1,6,5,2,4]$ to $[6,2,4,5,1,3]$ is
$$
[3,1,6,5,2,4]
\xrightarrow{t_{2,5}} [3,\textcolor{red}{2},6,5,\textcolor{red}{1},4]
\xrightarrow{t_{1,6}} [\textcolor{red}{4},2,6,5,1,\textcolor{red}{3}]
\xrightarrow{t_{1,3}} [\textcolor{red}{6},2,\textcolor{red}{4},5,1,3],
$$
Notice that $5,6,3$ are distinct numbers. 
\end{exa}

Then we use Theorem~\ref{T: Pieri} to deduce the following.
\begin{prop}
\label{P: Swap chains}
Pick $u, w \in S_n$, $k_1, k_2 \in [n-1]$
and $d_1, d_2 \in \mathbb{Z}_{\geq 0}$.
The number of chains from $u$ to $w$ compatible with $(k_1, k_2)$ and has weight $(d_1, d_2)$
matches the number of chains from $u$ to $w$ compatible with $(k_2, k_1)$ and has weight $(d_2, d_1)$.
\end{prop}
\begin{proof}
We may embed $S_n$ into $S_N$
for some $N$ large enough:
$v \in S_n$ is identified with
$v' \in S_N$ such that $v'(i) = v(i)$
if $i \in [n]$ and $v'(i) = i$ if $n < i \leq N$.
Now assume $u$ and $w$ are identified with
$u', w' \in S_N$.
Clearly, 
\begin{align*}
& |\{ (u, v, w): v \in S_n, \textrm{compatible with $(k_1, k_2)$, has weight $(d_1, d_2)$}\}|\\
=& |\{ (u', v', w'): v' \in S_N, \textrm{compatible with $(k_1, k_2)$, has weight $(d_1, d_2)$}\}|. 
\end{align*}

By Theorem~\ref{T: Pieri}, 
this
is the coefficient of $\fS_{w'}$ 
in $\fS_{u'} \times h_{d_1}(x_1, \cdots, x_{k_1}) \times h_{d_2}(x_1, \cdots, x_{k_2})$.
The proof is finished by the commutativity of 
polynomial multiplication.
\end{proof}

We are ready to introduce $(n-1)!$ 
many chain formulas for Schubert polynomials. 
\begin{cor}
\label{C: Bruhat formula}
Take $w \in S_n$ and $\gamma \in S_{n-1}$.
Then 
$\fS_w = \sum_{C} x^{\wt_\gamma(C)}$,
summing over all chains $C$ from $w$ 
to $w_0$ compatible with $\gamma$.
The summand $x^{\wt_\gamma(C)}$
is a monomial where the power of $x_i$
equals $n - i$ minus the length
of the increasing $i$-chain.
\end{cor}
\begin{proof}[Proof of Corollary~\ref{C: Bruhat formula}]
By~\ref{EQ: Bruhat chains}, 
we know Corollary~\ref{C: Bruhat formula} holds
when $\gamma = [1, \cdots, n-1]$.
Then the proof is finished by an induction
using Proposition~\ref{P: Swap chains}.
\end{proof}

We end this section by a question. 
In Proposition~\ref{P: Swap chains},
we observe there are two sets of the same size.
It would be natural to ask for an explicit bijection.
\begin{prob}
\label{prob: chain}
Find an explicit bijection
between the two set of chains in Proposition~\ref{P: Swap chains}. 
\end{prob}

In Corollary~\ref{C: growth 2}, 
we show Lenart's growth diagram~\cite{L_growth} 
solves Problem~\ref{prob: chain} in a special case.

\subsection{Fomin-Stanley construction}
\label{SS: FS}

Use $s_j$ to denote $t_{j, j+1}$.
A \definition{reduced word} of $w \in S_n$
is a word $i_1 i_2 \cdots i_l$ such that
$w = s_{i_1} \cdots s_{i_l}$ and $l$ is minimized.
One can read off a reduced word of $w$ 
from every $P \in \PD(w)$ as follows:
Go through its crossings from top to bottom
and right to left in each row.
For a crossing in row $r$ column $c$,
read off $r + c - 1$.
For instance, the PD in 
Example~\ref{E: PD} gives 
$41324$ which is a reduced word of $[2,5,1,4,3]$.

The \definition{nil-Coexter algebra} $\N_n$ is generated by $u_1, \cdots, u_{n-1}$
satisfying:
\begin{align*}
\begin{cases}
u_i^2 & = 0,\\
u_i u_j & =  u_j u_i \textrm{ if $|i-j| \geq 2$ },\\
u_i u_{i+1} u_i & =  u_{i+1} u_i u_{i+1} \textrm{ if $i \in [n-2]$.}
\end{cases}
\end{align*}

Consider $a = u_{i_1} \cdots u_{i_l} \in \N_n$, we have
$a \neq 0$ if and only if $i_1 \cdots i_l$
is a reduced word of some $w \in S_n$.
In this case,
$a = u_{j_1} \cdots u_{j_{l'}}$
if and only if $j_1 \cdots j_{l'}$
is a reduced word for the same $w$.
Fomin and Stanley~\cite{FS}
defined the following elements in $\Q[x_1, \cdots, x_{n-1}] \otimes \N_n$:
\begin{align*}
A_i(x_i) & := (1 + x_i u_{n-1})(1 + x_i u_{n-2}) \cdots (1 + x_i u_{i}) \textrm{ for $i \in [n-1]$, and}\\
\fS^\PD & := A_1(x_1) A_2(x_2) \cdots A_{n-1}(x_{n-1}).
\end{align*}
Combinatorially, 
$\fS^\PD = \sum_{P} x^{\wt(P)} u_{i_1}\cdots u_{i_l}$
where the sum runs over all PD
and $i_1 \cdots i_l$ is the reduced word read off from the PD.
Algebraically, Fomin and Stanley
showed that 
\begin{align}
\label{EQ: classical FS}
\fS^\PD = \sum_{w \in S_n} \fS_w u_{i_1}\cdots u_{i_l}, 
\end{align}
where $i_1 \cdots i_l$ 
is an arbitrary reduced word of $w$.
This formula would imply 
the PD formula.
Fomin and Stanley proved 
(\ref{EQ: classical FS})
by showing 
$\partial_i(\fS^\PD) = \fS^\PD u_i$
for any $i \in [n-1]$.
This equation then reduces
to 
$\partial_i(R_i(x_i) R_{i+1}(x_{i+1}))
= R_i(x_i) R_{i+1}(x_{i+1}) u_i$.
In \S\ref{S: FS},
we present the BPD analogue
of (\ref{EQ: classical FS})
and establish our equation
in a similar way.

\section{Bijetion between BPDs and chains}
\label{S: Chain}
In this section, we establish a bijection between 
$\BPD(w)$ and chains from $w$ to $w_0$
compatible with $(n-1, \cdots, 1)$.
Our description requires the following notions. 
\begin{defn}
\label{D: k se}
For $k = 0, 1,\cdots, n-1$, 
let $S_n^{k \searrow}$ be the subset of $S_n$
consisting of all $w$ such that 
$w(k+1) > w(k + 2) > \cdots > w(n)$.
Notice that a permutation $w \in S_n^{k \searrow}$
is completely determined by $w(1), \cdots w(k)$. 
\end{defn}

\begin{defn}
For $k = 0, \cdots, n-1$,
a \definition{$k$-row BPD} is a filling of 
a rectangular grid with $k$ rows and $n$ columns,
using tiles of BPDs. 
Similar to a BPD, pipes enter from the right edge
of each tile in column $n$ and cannot cross more than once.
A $k$-row BPD is associated with a permutation 
$w \in S_n^{k \searrow}$.
For each $i \in [k]$, we let $w(i)$
be the column index of the tile in row $k$
where the pipe from row $i$ exits from the bottom.
\end{defn}

\begin{exa}
Suppose $n = 6$.
The following is a $3$-row BPD associated with
$w \in S_6^{3 \searrow}$ with one-line notation $316542$.
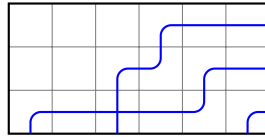
\[ 
\begin{tikzpicture}[x=1.5em,y=1.5em,thick,rounded corners, color = blue]
\draw[step=1,gray,thin] (0,0) grid (6,3);
\draw[color=black, thick, sharp corners] (0,0) rectangle (6,3);
\draw(0.5, 0)--(0.5,0.5)--(4.5,0.5)--(4.5,1.5)--(6,1.5);
\draw(2.5, 0)--(2.5,1.5)--(3.5,1.5)--(3.5,2.5)--(6,2.5);
\draw(5.5, 0)--(5.5,0.5)--(6,0.5);
\end{tikzpicture}
\]
\end{exa}

Intuitively, a $k$-row BPD is just the 
first $k$ rows of a BPD. 
Now we can introduce the bijection.

\begin{defn}
Take $w \in S_n$ and $D \in \BPD(w)$.
For $k = 0, \cdots, n-1$,
let $u_k$ be the permutation
of the $k$-row BPD obtained from
keeping the first $k$-rows of $D$.
Define $\chain_\BPD(D)$ as the chain
$(u_{n-1}, \cdots, u_0)$.
\end{defn}

\begin{exa}
\label{E: BPD to chain}
Assume $n = 6$.
Take $D \in \BPD([2,1,6,5,3,4])$
as depicted below. 
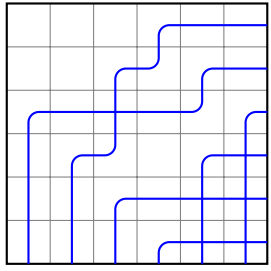
\[
\begin{tikzpicture}[x=1.5em,y=1.5em,thick,rounded corners, color = blue]
\draw[step=1,gray,thin] (0,0) grid (6,6);
\draw[color=black, thick, sharp corners] (0,0) rectangle (6,6);
\draw(0.5, 0)--(0.5,3.5)--(4.5,3.5)--(4.5,4.5)--(6,4.5);
\draw(1.5, 0)--(1.5,2.5)--(2.5,2.5)--(2.5,4.5)--(3.5,4.5)--(3.5,5.5)--(6,5.5);
\draw(2.5, 0)--(2.5,1.5)--(6,1.5);
\draw(3.5, 0)--(3.5,0.5)--(6,0.5);
\draw(4.5, 0)--(4.5,2.5)--(6,2.5);
\draw(5.5, 0)--(5.5,3.5)--(6,3.5);
\end{tikzpicture}
\]
Then $\chain_\BPD(D)$ is the following chain
$$
([2,1,6,5,3,4], [2,1,6,5,4,3],[3,1,6,5,4,2],
[3,5,6,4,2,1], [4,6,5,3,2,1], [6,5,4,3,2,1]).
$$
\end{exa}

The main theorem of this section is an analogue of 
Theorem~\ref{T: PD to chain}.
It shows our map $\chain_\BPD$
parallels the map $\chain_\PD$ 
of Lenart and Sottile.

\begin{thm}
\label{T: BPD to chain}
Take $w \in S_n$.
The map $\chain_\BPD$
is a bijection from $\BPD(w)$
to chains from $w$ to $w_0$
that is compatible with $(n-1, \cdots, 1)$.
Moreover, if $\chain_\BPD(D) = C$,
then $\wt(C)$ is the reverse 
of $(n-1, \cdots, 1) - \wt(D)$.
\end{thm}

The key step to prove it is to explain why $\chain_\BPD$
yields a chain compatible with $(n-1, \cdots, 1)$
and construct its inverse.
We start by characterizing what permutations
can show up in a chain from $w$ to $w_0$
compatible with $(n-1, \cdots, 1)$.

\begin{lemma}
\label{L: Stay in S^k (perm)}
Take $k \in [n-1]$
and $w \in S_n^{k \searrow}$.
Take $k' \leq k$.
For any $u$ such that $u \leq_{k'} w$,
we must have $u \in S_n^{k \searrow}$.
\end{lemma}
\begin{proof}
It is enough to assume $u = w t_{i,j} \lessdot_k w$
for some $i \leq k' < j$.
If $j \leq k$, we are done. 
Otherwise, $i \leq k < j$.
To check $u \in S_n^{k \searrow}$, 
we just need to make sure the following two statements:
\begin{itemize}
\item If $j < n$, we need to make sure 
$u(j) > u(j+1)$:
Since $w t_{i, j}  \lessdot w$,
$w(i) > w(j)$.
Since $w \in S_n^{k \searrow}$, 
$w(j) > w(j+1)$.
Thus, $u(j) = w(i) > w(j) > w(j+1) = u(j+1)$.
\item If $j-1 > k$,
we need to make sure $u(j-1) > u(j)$:
By $w \in S_n^{k \searrow}$,
we have $w(j-1) > w(j)$.
Then by $w t_{i, j}  \lessdot w$,
we have $w(j-1) > w(i)$.
Thus, 
$u(j-1) = w(j-1) > w(i) = u(j)$.\qedhere
\end{itemize}
\end{proof}

\begin{cor}
\label{C: chain in se}
Let $(u_{n-1}, \cdots, u_0)$ be a chain
compatible with $(n-1, \cdots, 1)$
and $u_0 = w_0$.
Then $u_i \in S_n^{i\searrow}$.
\end{cor}
\begin{proof}
Prove by induction on $i$.
First, $u_0 = w_0 \in S_n^{0\searrow}$.
Then for $i \in [n-1]$,
there is an $i$-chain from $u_i$
to $u_{i-1}$,
which is in $S_n^{i-1 \searrow}$
by our inductive hypothesis. 
Thus, $u_{i-1} \in S_n^{i \searrow}$.
Then $u_i \in  S_n^{i \searrow}$
by Lemma~\ref{L: Stay in S^k (perm)}.
\end{proof}

We also need the following observation
that recasts a condition in
the definition of $k$-row BPDs.
\begin{rem}
\label{R: BPD double crossing}
Notice that the requirement
``two pipes cannot cross more than once''
in the definition of BPD and $k$-row BPD
is equivalent to the following:
For each $\ptile$,
if the horizontal (resp. vertical) pipe
is from row $i$ (resp. $j$),
then $i > j$.
\end{rem}

Next, we provide an explicit algorithm to construct the increasing
$k$-chain from $u_k$ to $u_{k-1}$
when $\chain_\BPD(D) = (u_{n-1}, \cdots, u_1, u_0)$.

\begin{lemma}
\label{l: One row to chain}
Take $k \in [n-1]$ and let $D$ be a $k$-row BPD
associated with $u$.
Suppose the first $(k-1)$ rows of $D$ 
is a $(k-1)$-row BPD associated with $w$.
Let $c_1< \cdots< c_m$ be the 
column indices of $\rtile$ and $\htile$
in row $k$ of $D$.
Suppose the pipe in $(k,c_i)$ is from row $r_i$.
Then
\begin{equation}
\label{eq: chain from a row}
u \xrightarrow{t_{r_1, w^{-1}(c_1)}} \cdots \xrightarrow{t_{r_{m-1}, w^{-1}(c_{m-1})}} w
\end{equation}
is an increasing $k$-chain from $u$ to $w$.
\end{lemma}

In other words, 
$u$ is obtained from $w$
by exchanging the value $c_i$ 
with the $r_i\textsuperscript{th}$ number
for $i = m-1, \cdots, 1$.

\begin{exa}
The $1$-row BPD obtained by keeping 
the first row of $D$ in Example~\ref{E: BPD to chain} is associated with $[4,6,5,3,2,1]$.
The $2$-row BPD obtained by keeping 
the first two rows 
is associated with $[3,5,6,4,2,1]$.
Lemma~\ref{l: One row to chain} asserts
that there is an increasing $2$-chain
from $[3,5,6,4,2,1]$ to $[4,6,5,3,2,1]$.
To find this chain, 
we look at row $2$ of $D$.
The $\rtile$ and $\htile$ are in column 
$3, 5$ and $6$.
The tile $(2,3)$ contains the pipe from row $1$
and the tile $(2,4)$ contains the pipe from 
row $2$.
By~\eqref{eq: chain from a row},
$$
[3,5,6,4,2,1] \xrightarrow{t_{1,4}}
[4,5,6,3,2,1] \xrightarrow{t_{2,3}}
[4,6,5,3,2,1]
$$
is the increasing $2$-chain. 
\end{exa}

\begin{proof}
Define $w_m = w$
and $w_i = w_{i+1} t_{r_i, w^{-1}(c_i)}$.
We need to make sure $w_1 = u$ 
and 
\begin{equation}
w_1 \xrightarrow{t_{r_1, w^{-1}(c_1)}} w_2 \xrightarrow{t_{r_2, w^{-1}(c_2)}}
\cdots w_{m-1} 
\xrightarrow{t_{r_{m-1}, w^{-1}(c_{m-1})}} w_m
\end{equation}
is an increasing $k$-chain.
We start by making some useful observations.
\begin{itemize}
\item Observation 1: The tile
$(k, c_m)$ is the rightmost $\rtile$ or $\htile$
in row $k$,
so it must contain the pipe from row $k$.
\item Observation 2: The tile $(k, c_m)$ is the rightmost
tile in row $k$ without a pipe entering from the top.
By $w \in S_n^{k-1 \searrow}$, $w(k) = c_m$.
\item Observation 3: For $i \in [m-1]$, the value $c_i$
is not among the first $k$ numbers of $w_{i+1}$:
Because $(k, c_i)$ has no pipe
entering from above,
we know $c_i$ is not among the first
$k-1$ numbers of $w$. 
By Observation 1, $w(k) = c_m$,
so $c_i$ is not among the first $k$ numbers of $w = w_m$.
Then $w_{i+1}$ is obtained from $w$ by 
repetitively swapping someone among the first $k$
numbers and someone larger than $c_i$,
so $c_i$ stays at the same place in $w_{i+1}$ and $w$.
\item Observation 4: If the pipe from row $r$ goes straight down
in row $k$ without making a turn, 
then $w_i(r)$ stays the same for all $i$.
\end{itemize}

For each $c_i$ with $i \in [m-1]$, 
we can find smallest $d_i > c_i$
such that tiles between $(k, c_i)$
and $(k,d_i)$ are $\ptile$.
Clearly, 
$(k,d_i)$ is $\htile$ or $\jtile$,
containing the same pipe as $(k, c_i)$.
We claim $w_{i+1}(r_i) = d_i$.
This claim can be checked through a case study:
\begin{itemize}
\item If $(k, d_i)$ is a $\jtile$,
then the pipe from row $r_i$
enters row $k$ from $(k, d_i)$,
so $w(r_i) = d_i$.
Since the pipe from row $r_i$ 
does not appear on the right of $(k, d_i)$,
we know $w_{i+1}(r_i) = w(r_i) = d_i$.
\item Otherwise, $(k, d_{i})$ is a $\htile$,
so $d_{i} = c_{i+1}$.
If $i+1 < m$,
we know $w_{i+1}$ 
is obtained from $w_{i+2}$
by swapping the value $c_{i+1}$
to the $i\textsuperscript{th}$ entry,
so $w_{i+1}(r_i) = c_{i+1}$.
If $i+1 = m$,
by Observation~$1$, $r_{i} = k$.
Then $w_m(r_i) = w(k) = c_m$
by Observation~$2$.
\end{itemize}

Thus, we know $w_i$ is obtained from
$w_{i+1}$ by exchanging the values
$c_i$ and $d_i$.
Next, we show this is a saturated $k$-chain
by checking $w_i \lessdot_k w_{i+1}$.
When we swap $c_i$ and $d_i$ in $w_{i+1}$,
we know $d_i$ is at the $r_i\textsuperscript{th}$
position with $r_i \leq k$.
By Observation 3,
$c_i$ is not among the first $k$ numbers
in $w_{i+1}$.
It remains to check for any $c_i < j < d_i$,
$j$ is on the left of $d_i$ in $w_{i+1}$.
By our assumption, $(k, j)$ is $\ptile$.
In this tile, the horizontal
pipe is from row $r_i$ and the vertical pipe
is from row $w^{-1}(j)$.
By Remark~\ref{R: BPD double crossing},
we have $w^{-1}(j) < r_i$.
By Observation 4, 
$w_{i+1}^{-1}(j) < r_i$.
In other words, 
$j$ appears on the left of $d_i$ in $w_{i+1}$.
We have $w_i \lessdot_k w_{i+1}$.

To see the chain is increasing, 
just notice that the smaller number 
being swapped are $c_1, \cdots, c_{d-1}$,
which are increasing.
Finally, we verify $u = w_1$.
We know $u \in S_n^{k \searrow}$ by definition.
Also, we have shown $w_1 \leq_k w_m$,
so Lemma~\ref{L: Stay in S^k (perm)} implies
$w_1 \in S_n^{k \searrow}$.
It remains to check $u(i) = w_1(i)$ 
for all $i \in [k]$.
If the pipe from row $i$ goes straight down
in row $k$, then $w_1(i) = w(i)$ by observation $4$.
Correspondingly, we have $u(i) = w(i)$.
Otherwise, the pipe from row $i$ 
has a $\rtile$ in row $k$,
say at column $c$.
Then $w_1(i) = c$.
Correspondingly, we have $u(i) = c$.
\end{proof}

We may further simplify 
the description of the chain in~\eqref{eq: chain from a row}
as follows. 
This simplification will be used in the next section. 

\begin{cor}
\label{c: One row to chain}
Take $k \in [n-1]$ and let $D$ be a $k$-row BPD
associated with $u$.
Suppose the first $(k-1)$ rows of $D$ 
is a $(k-1)$-row BPD associated with $w$.
Let $c_1< \cdots< c_m$ be the 
column indices of $\rtile$ and $\htile$
in row $k$ of $D$.
Suppose $(k,c_i)$ is the $j_i\textsuperscript{th}$
leftmost tile in row $k$ that is 
$\rtile$, $\htile$ or $\btile$.
Then
$$
u \xrightarrow{t_{r_1, n +1 - j_1}} \cdots \xrightarrow{t_{r_{m-1}, n + 1 - j_{m-1}}} w
$$
is an increasing $k$-chain from $u$ to $w$.    
\end{cor}
\begin{proof}
By Lemma~\ref{l: One row to chain},
we only need to check $w^{-1}(c_i) = n + 1 - j_{i}$
for each $i \in [m-1]$.
Notice that if a tile $(k, c)$ is $\rtile$, 
$\htile$, or $\btile$,
then it has no pipe entering from its top, 
so $c$ is not among the first $k-1$ numbers in $w$.
By definition, $(k, c_i)$ is such a tile 
and there are $j_i - 1$ such tiles on its left. 
In other words, 
$c_i$ is the $j_i\textsuperscript{th}$
smallest number among $w(k), \cdots, w(n)$.
Finally, by $w \in S_n^{k-1 \searrow}$,
we know $w(n+1-j_i) = c_i$.
\end{proof}

Following this corollary, 
we make an observation that parallels
Remark~\ref{R: Double PD} for PDs. 
\begin{rem}
\label{R: Double BPD}
Suppose $\chain_\BPD(B) = (u_{n-1}, \cdots, u_0)$.
Then by Lemma~\ref{l: One row to chain},
there is a increasing $k$-chain
from $u_{k}$ to $u_{k-1}$.
We know a number $j$ appears among the last
$n - k + 1$ numbers of $u_{k-1}$
if the tile $(k, j)$ is $\rtile$, $\htile$
or $\btile$.
The column index of the rightmost such tile
in row $k$ is $u_{k-1}(k)$.
By Corollary~\ref{c: One row to chain}, 
the tile $(k, j)$ is a $\btile$
if and only if it is among
the last $n - k$ numbers of 
$u_{k-1}$ and remains unswapped
in the increasing $k$-chain from $u_{k}$
to $u_{k-1}$.
\end{rem}

Now we describe the inverse of $\chain_\BPD$.
We start with the following lemma.
\begin{lemma}
\label{l: Chain restriction}
Take $k \in [n-1]$ and
$w \in S_{n}^{k \searrow}$.
Assume there is an increasing 
$k$-chain from some $u$ to $w$.
Suppose the chain swaps the values $a < b$,
and assume after the swap, 
$b$ is at index $i$.
We make three claims. 
\begin{itemize}
\item The number $b$ is not among
the first $k$ numbers of $u$.
\item For any $c$ such that $a< c < b$,
we have $u^{-1}(c) < i$.
\item For any $c$ such that $a< c < b$,
the value $c$ stays unswapped
throughout the chain.
\end{itemize}
\end{lemma}
\begin{proof}
Let $u' \lessdot w'$
be the two permutations in the chain
where $w'$ is obtained from $u'$
by swapping values $a$ and $b$.
We prove the claims separately. 
\begin{itemize}
\item Suppose $b$ is among the first $k$
numbers in $u$.
We know $b$ is not among the first $k$
numbers in $u'$,
so it must be swapped with some 
larger number in the chain segment
from $u$ to $u'$.
This contradicts to the chain being increasing. 
\item 
Since $a < c < b$, 
in $w'$, 
the value $c$ is either on the 
left of $b$ or on the right of $a$.
We show the latter case is impossible.
Since $w' \leq_k w$,
Lemma~\ref{L: Stay in S^k (perm)}
implies $w' \in S_n^{k \searrow}$.
We know $a$ is not among the 
first $k$ numbers of $w'$,
then $c$ cannot appear on its right since $a < c$.

Now we know there is $j < i$ such that $u'(j) = w'(j) = c$.
Since the chain is increasing, 
$u$ is obtained from $u'$ 
by doing transpositions where the smaller number being 
swapped is less than $a$.
In the chain segment from $u$ to $u'$,
the value $a$ stays at index $i$, 
so $c$ cannot be swapped with any number smaller than $a$.
Thus, $u(j) = u'(j) = c$.

\item 
From the second claim, we know $c$ is among 
the first $k$ numbers in $u$.
Suppose $c$ is swapped with $d$
in the chain.
By the first claim, we have $c < d$
where $d$ is not among the first 
$k$ numbers in $u$.
If $b < d$, we have $c < b < d$.
Then by our second claim,
$b$ is among the first $k$ numbers in $u$,
which contradicts to our first claim.
Otherwise, 
$a < d < b$.
Our second claim implies
$d$ is among the first $k$ numbers of $u$,
which is also a contradiction.
\qedhere
\end{itemize}
\end{proof}

Now we provide an algorithm that 
constructs a row of BPD. 
\begin{defn}
\label{D: one row BPD}
Take $k \in [n-1]$ and
$w \in S_{n}^{k \searrow}$.
Assume there is an increasing 
$k$-chain from some $u$ to $w$.
We construct a row of $n$ tiles
using the tiles of a BPD. 
The tile $(1,c)$
can be determined as follows. 
First, we have three cases if 
the chain swaps the value $c$. 
\begin{itemize}
\item If the chain swaps $c$
with someone larger than it,
but not with someone smaller than it, 
this tile is $\rtile$.
\item If the chain swaps $c$
both with someone larger than it and
someone smaller than it, 
this tile is $\htile$.
\item If the chain swaps $c$
with someone smaller than it,
but not with someone larger than it, 
this tile is $\rtile$.
\end{itemize}
If $c$ stays unswapped in the chain,
we also have three cases. 
\begin{itemize}
\item If $c$ is not among the 
first $k$ numbers of $w$,
this tile is $\btile$.
\item Otherwise, if the chain 
ever swaps values $a, b$
with $a < c < b$,
then the tile is $\ptile$.
\item Otherwise, the tile is $\vtile$.
\end{itemize}
\end{defn}

\begin{exa}
Let $n = 6$ and $k = 2$. 
Consider the following increasing 
$2$-chain:
$$
[3,5,6,4,2,1] \xrightarrow{t_{1,4}}
[4,5,6,3,2,1] \xrightarrow{t_{2,3}}
[4,6,5,3,2,1].
$$
Since $w = [4,6,5,3,2,1] \in S_6^{2 \searrow}$,
we may construct a row of tilings
using Definition~\ref{D: one row BPD}.
Since $1$ and $2$ are unswapped and not
among the first $2$ numbers of $w$,
the first two tiles are $\btile$.
Since $3$ is only swapped with $4$,
the next tile is $\rtile$.
Since $4$ is only swapped with $3$,
the next tile is $\jtile$.
Since $5$ is only swapped with $6$,
the next tile is $\rtile$.
Since $6$ is only swapped with $5$,
the next tile is $\jtile$.
Thus, we obtain:
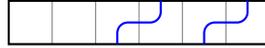
\[ 
\begin{tikzpicture}[x=1.5em,y=1.5em,thick,rounded corners, color = blue]
\draw[step=1,gray,thin] (0,0) grid (6,1);
\draw[color=black, thick, sharp corners] (0,0) rectangle (6,1);
\draw(2.5, 0)--(2.5,0.5)--(3.5,0.5)--(3.5,1);
\draw(4.5, 0)--(4.5,0.5)--(5.5,0.5)--(5.5,1);
\end{tikzpicture}
\]

As another example, 
take $k = 3$
and consider the increasing $3$-chain:
$$
[3,6,1,5,4,2] \xrightarrow{t_{3,6}}
[3,6,2,5,4,1] \xrightarrow{t_{3,5}}
[3,6,4,5,2,1] \xrightarrow{t_{3,4}}
[3,6,5,4,2,1].
$$
Since $[3,6,5,4,2,1] \in S_6^{3 \searrow}$,
we may construct a row of tilings
using Definition~\ref{D: one row BPD}:
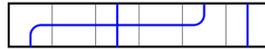
\[ 
\begin{tikzpicture}[x=1.5em,y=1.5em,thick,rounded corners, color = blue]
\draw[step=1,gray,thin] (0,0) grid (6,1);
\draw[color=black, thick, sharp corners] (0,0) rectangle (6,1);
\draw(.5, 0)--(.5,0.5)--(4.5,0.5)--(4.5,1);
\draw(2.5, 0)--(2.5, 1);
\draw(5.5, 0)--(5.5, 1);
\end{tikzpicture}
\]
\end{exa}

\begin{lemma}
\label{L: one row BPD}
Let $k ,u$ and $w$ 
be defined in Definition~\ref{D: one row BPD}.
Then the row of tiles we construct
satisfies the following:
\begin{itemize}
\item The pipes are connected. 
\item For each $i \in [k]$,
there is a pipe entering the row
from the top of $(1, w(i))$ and exiting from
the bottom of $(1, u(i))$.
\item For any $\ptile$,
assume the horizontal pipe
enters the row from $(1, w(i))$
and the vertical pipe enters the row from $(1,w(j))$,
then $j < i$.
\end{itemize}
\end{lemma}
\begin{proof}
\noindent
\begin{itemize}
\item Consider the tile $(1,c)$ with $c \in [n-1]$.
We need to show that $(1,c)$ has a pipe connecting 
to the right 
if and only if $(1,c+1)$ has a pipe connecting 
to the left. 
We check the forward direction and the other direction
can be checked similarly. 

Suppose $(1,c)$ is $\rtile$ or $\htile$,
then we know the chain swaps the value $c$ 
with some $c' > c$.
If $c' = c+1$, 
then $(1, c+1)$ is $\jtile$ or $\htile$,
both are connected to the left.
If $c' > c+1$,
by Lemma~\ref{l: Chain restriction}, 
we know $c+1$ is among the first $k$
numbers of $w$ and remains unswapped in the chain.
Thus, $(1, c+1)$ is $\ptile$.

Finally, suppose $(1,c)$ is $\ptile$,
we know the chain swaps $a, b$ 
such that $a < c < b$.
If $b = c+1$, 
then $(1, c+1)$ is $\jtile$ or $\htile$,
both are connected to the left.
If $b > c+1$,
by Lemma~\ref{l: Chain restriction}, 
we know $c+1$ is among the first $k$
numbers of $w$ and remains unswapped in the chain.
Thus, $(1, c+1)$ is $\ptile$.

\item 
Fix $i \in [k]$.
If $w(i)$ is not swapped by the chain,
we know $(1, w(i))$ is $\ptile$ or $\vtile$
and our claim is immediate. 
Otherwise, we may find a sequence of numbers 
$c_1 < \cdots < c_m = w(i)$
such that $c_j$ is swapped with $c_{j+1}$
and $c_1$ is not swapped with any smaller number.
Thus, $u(i) = c_1$.
By our construction of the row,
the tile $(1,c_1)$ is $\rtile$,
the tile $(1,c_m)$ is $\jtile$,
and tiles between them are $\htile$ or $\ptile$.
Thus, there is a pipe going from 
the top of $(1,w(i))$ to 
the bottom of $(1, u(i))$.
\item 
Let $(1, c)$ be a $\ptile$
where the horizontal pipe 
enters the row from $(1, w(i))$.
We may find a sequence of numbers 
$c_1 < \cdots < c_m = w(i)$
such that $c_j$ is swapped with $c_{j+1}$
and $c_1$ is not swapped with any smaller number.
We know $c_t < c < c_{t+1}$ for some $t$.
In the chain,
after it swaps $c_t$ and $c_{t+1}$,
we know $c_{t+1}$ is at index $i$.
By Lemma~\ref{l: Chain restriction},
$w^{-1}(c) < i$. \qedhere
\end{itemize}
\end{proof}

\begin{lemma}
\label{L: one row BPD2}
Let $D$ be a $(k-1)$-row BPD 
associated with $w \in S_n^{k-1 \searrow}$.
Assume there is an increasing $k$-chain from 
some permutation $u$ to $w$.
Then there exists a $k$-row BPD
associated with $u$ such that its
first $(k-1)$ rows agree with $D$.
\end{lemma}
\begin{proof}
Notice that 
$w \in S_n^{k-1 \searrow} \subset S_n^{k \searrow}$.
Consider the row of tilings from
Definition~\ref{D: one row BPD}.
By the second claim in Lemma~\ref{L: one row BPD},
the tile $(1, w(i))$ has a pipe entering from the top
for each $i \in [k]$.
Then $(1, w(k))$ cannot be $\ptile$:
If so, the the horizontal pipe enters the 
row from $(1, w(i))$ for some $i < k$,
contradicting the last claim
in Lemma~\ref{L: one row BPD}.

Next, by $w \in S_n^{k-1 \searrow}$,
numbers larger than $w(k)$ 
are among the first $k-1$ numbers in $w$.
Thus, there is a pipe entering from the top
in each $(1, j)$ with $w(k) < j \leq n$.
Since $(1, w(k))$ is not a $\ptile$,
all tiles on its right must be $\vtile$.

Finally, we change $(1, w(k))$
from $\jtile$ to $\htile$
or from $\vtile$ to $\rtile$.
Then we change $(1, j)$ from
$\vtile$ to $\ptile$
for each $j > w(k)$.
We append this row
on the bottom of $D$ and 
let $D'$ be the resulting tiling. 
We check $D'$ is a $k$-row BPD associated with $u$.
By the second claim in Lemma~\ref{L: one row BPD}, 
the pipe from row $i$ in $D'$ ends at column $u(i)$.
If $(k, c)$ is a $\ptile$ in $D'$,
say the vertical pipe comes from row $j$
and the horizontal pipe comes from row $i$.
By the third claim in Lemma~\ref{L: one row BPD}, 
$j < i$.
\end{proof}

We are prepared to prove Theorem~\ref{T: BPD to chain}.
\begin{proof}[Proof of Theorem~\ref{T: BPD to chain}]
We first show $\chain_\BPD$
is a well-defined map.
Take $D \in \BPD(w)$.
Let $D_k$ be the $k$-row BPD obtained by keeping 
the first $k$ rows of $D$.
Let $u_k$ be the permutation associated with 
$D_k$.
Clearly, 
$u_0 = w_0$ and $u_{n-1} = w$.
By Lemma~\ref{l: One row to chain},
there is an increasing 
$k$ chain from $u_k$ 
to $u_{k-1}$ for $k \in [n-1]$.
Thus, 
$\chain_\BPD(D) = 
(u_{n-1}, \cdots, u_0)$
is compatible with $(n-1, \cdots, 1)$.

Next, we show $\chain_\BPD$
respects the weight. 
For $k \in [n-1]$,
consider the $k\textsuperscript{th}$
entry of $\wt(D)$,
which is the number of $\btile$ in row $k$ of $D$.
By Remark~\ref{R: Double BPD},
it equals the number of values that 
are among the last $n-k$ numbers in
$u_{k-1}$ and remains unswapped 
in the increasing chain from $u_k$ to $u_{k-1}$.
By Remark~\ref{R: wt of inc},
this equals $n - k$ minus the length
of the increasing chain.
In other words, 
it is the $k\textsuperscript{th}$ rightmost
entry in $(1, \cdots, n-1) - \wt(\chain_\BPD(D))$.

Then we show $\chain_\BPD$ is injective.
Suppose $\chain_\BPD(D') = \chain_\BPD(D)$.
Then for each $k \in [n]$,
the starting points and destinations of each pipe
are the same in row $k$ of $D$ and $D'$.
Clearly, row $k$ of $D$ and $D'$ must agree, 
so $D = D'$.

Finally, we construct the inverse of $\chain_\BPD$.
Let $(u_{n-1}, \cdots, u_{0})$ be a chain 
from $w$ to $w_0$ compatible with $(n-1, \cdots, 1)$.
By Corollary~\ref{C: chain in se},
we know $u_k \in S_n^{k \searrow}$.
We construct $D_0, \cdots, D_{n-1}$
where $D_k$ is a $k$-row BPD associated with $u_k$. 
We start by setting $D_0$ as the only $0$-row BPD,
which has no tiling at all. 
Then we may find a $D_{k+1}$ from $D_k$
using Lemma~\ref{L: one row BPD2}.
After we get $D_{n-1}$,
there is clearly a way to append 
a row on the bottom and obtain a BPD $D \in \BPD(u_{n-1})$.
The first $k$ rows of $D$ agrees with $D_k$,
so $\chain_\BPD(D) = (u_{n-1}, \cdots, u_0)$.
Thus, this map $(u_{n-1}, \cdots, u_0) \mapsto D$ is a right inverse of $\chain_\BPD(\cdot)$.
Since $\chain_\BPD(\cdot)$ is injective,
this map is its inverse. 
\end{proof}

\section{Analogue of Fomin-Stanley construction on BPD}
\label{S: FS}

Our first application of Theorem~\ref{T: BPD to chain}
is a BPD analogue of Fomin-Stanley construction. 
Instead of the nil-Coexter algebra $\N_n$, our construction
uses the \definition{Fomin-Kirillov algebra}~\cite{FK-algebra}.
In Section~\ref{SS: FK}, we cover some
background regarding the Fomin-Kirillov algebra.
Then in Section~\ref{SS: Flagged tableaux},
we introduce a new combinatorial object which we call
flagged tableaux.
They naturally encode the information of a BPD.
We define $\FT(w)$ as a set  
of certain flagged tableaux and build a $wt(\cdot)$-preserving 
bijection from $\BPD(w)$ to $\FT(w)$.
In Section~\ref{SS: construction},
we construct $\fS^{\BPD} \in \mathbb{Q}[x_1, \cdots, x_{n-1}][S_n]$.
By the following Proposition, 
$\fS^\BPD$ is the generating function of
flagged tableaux, or equivalently BPDs.
\begin{prop}
\label{P: generating function}
We have
$$
\fS^\BPD = \sum_{w \in S_n} \sum_{T \in \FT(w)} x^{\wt(T)}w
= \sum_{w \in S_n} \sum_{D \in \BPD(w)} x^{\wt(D)}w.
$$
\end{prop}
In Section~\ref{SS: FS Proofs},
we present an operator theoretic proof 
of the following.
\begin{thm}
\label{T: FS on BPD}
We have 
$\fS^\BPD = \sum_{w \in S_n} \fS_{w} w$.  
\end{thm}
Thus, we obtain an alternative proof
of the BPD formula $\fS_w = \sum_{D \in \BPD(w)} x^{\wt(D)}$.
Our proof heavily relies on algebraic properties of
certain elements in $\N_n$
known as the Dunkl elements. 
Finally, in Section~\ref{SS: Commuative},
we write $\fS^\BPD$ as a product of terms that 
are commutative with each other.
It turns out that this formulation is related to
the Cauchy's identity of Schubert polynomials. 

\subsection{Fomin-Kirillov algebra}
\label{SS: FK}
The Fomin-Kirillov algebra $\E_n$ is
generated by $\{d_{i,j}: 1 \leq i < j \leq n\}$ satisfying: 
\begin{align*}
\begin{cases}
d_{i,j}^2 & = 0 \textrm{ if $i < j$ }, \\
d_{i,j} d_{j,k}  & = d_{i,k} d_{i,j} + d_{j,k} d_{i,k}
\textrm{ if $i < j < k$ },\\
d_{j,k} d_{i,j}  & = d_{i,j} d_{i,k} + d_{i,k} d_{j,k}
\textrm{ if $i < j < k$ },\\
d_{i,j} d_{k,l}  & = d_{k,l} d_{i,j} 
\quad \textrm{if $i < j$, $k < l$ and $i, j, k, l$ distinct.}\\    
\end{cases}
\end{align*}

Fomin and Kirillov~\cite{FK-algebra} described an action of $\E_n$ on $\Q[S_n]$.
In this paper, we adopt a slightly different convention
and consider a right action of $\E_n$ on $\Q[S_n]$.
For $w \in S_n$,
$$
w \odot d_{i,j} := \begin{cases}
w t_{i,j} & \text{if $wt_{i,j} \lessdot w$} \\
0 & \text{otherwise.}
\end{cases}
$$

Fomin and Kirillov~\cite{FK-algebra} defined the \definition{Dunkl element}
$\theta_i := -\sum_{j < i} d_{j, i} + \sum_{j > i} d_{i,j} \in \E_n$
for $i \in [n]$.
\begin{cor}[{\cite[Corollary. 5.2]{FK-algebra}}]
The Dunkl elements $\theta_1, \cdots, \theta_n$
commute with each other.    
\end{cor}

One way to understand the Dunkl element $\theta_i$ 
is that it models the multiplication of $x_i$
with a Schubert polynomial.
Recall a well-known reformulation of the Monk's rule~\eqref{eq: Monk}:
For $w \in S_n$ and $i \in [n-1]$,
if there is $j>i$ with $w(j) > w(i)$,
then
\begin{equation}
\label{eq: Monk2}
\fS_w x_i = - \sum_{j < i: \: w \lessdot wt_{j,i}} \fS_{wt_{j,i}}
+ \sum_{j > i: \: w \lessdot wt_{i,j}} \fS_{wt_{i,j}}.  \end{equation}
In other words, 
for such $w$ and $i$,
the coefficient 
of $\fS_u$ in $\fS_w x_i$
agrees with the coefficient of $u w_0$
in $-w w_0 \odot \theta_{n-i + 1}$.
A special case of~\eqref{eq: Monk2} is the 
famous ``transition formula''. 
Fix any $w \in S_n$ that is not $[1, \cdots, n]$.
Let $i$ be the largest such that 
there exists $j > i$ with $w(j) < w(i)$.
Then find the largest such $j$.
Following~\eqref{eq: Monk2}, 
one has
\begin{equation}
\label{eq: transition}
\fS_w = \fS_{w t_{i,j}} x_i
+ \sum_{h < i: \: w t_{i,j} \lessdot wt_{i,j}t_{h,i}} \fS_{wt_{i,j}t_{h,i}}.  
\end{equation}
Using the transition formula, 
we derive a lemma that will be useful later.

\begin{lemma}
\label{L: Dunkl Monk}
For any $w \in S_n$,
we may let $\fS_w(-\theta_n, \cdots, -\theta_2)$ 
be an element in $\N_n$ obtained by 
replacing $x_i$ in $\fS_w$ by $-\theta_{n + 1 - i}$.
Then $$w_0 \odot \fS_w(-\theta_n, \cdots, -\theta_2)
= w w_0.$$
\end{lemma}
\begin{proof}
Prove by induction on $\ell(w)$.
The base case, when $w = [1, \cdots, n]$, is immediate.
Now for $w \in S_n$ that is not $[1, \cdots, n]$,
we let $i, j$ be the numbers in~\eqref{eq: transition}.
We have
\begin{align*}
& w_0 \odot \fS_w(-\theta_n, \cdots, -\theta_2)\\
= & w_0 \odot -\fS_{w t_{i,j}}(-\theta_n, \cdots, -\theta_2) \theta_{n - i + 1}
+ \sum_{h < i: \: w t_{i,j} \lessdot wt_{i,j}t_{h,i}} w_0 \odot \fS_{wt_{i,j}t_{h,i}}(-\theta_n, \cdots, -\theta_2) \\
= & - w t_{i,j} w_0 \odot \theta_{n - i + 1}
+ \sum_{h < i: \: w t_{i,j} \lessdot wt_{i,j}t_{h,i}} wt_{i,j}t_{h,i} w_0 = w w_0 \qedhere
\end{align*}
\end{proof}

\subsection{Encoding BPDs as flagged tableaux}
\label{SS: Flagged tableaux}
We will encode each 
BPD as the following combinatorial object.
\begin{defn}
A \definition{flagged tableau}
is a staircase grid
with a cell in row $i$
column $j$ if $i + j \leq n$.
Moreover, each cell in row $i$ is empty
or filled with a number in $[i]$.    
\end{defn}

We define an encoding map $\Phi$
from $\BPD(w)$ to the set of flagged tableaux.

\begin{defn}
Take $D \in \BPD(w)$ for some $w \in S_n$.
To fill the cell in row $i$ column 
$j$ of $\Phi(D)$,
we look at the 
$j^\textsuperscript{th}$
leftmost $\btile$, 
$\htile$ or $\rtile$
in row $i$ of $D$.
If it is a blank,
we leave the cell in 
$\Phi(D)$ unfilled. 
Otherwise, it contains a pipe from row $r$ for some $r \leq i$.
We fill the cell
in $\Phi(D)$ by $r$.
\end{defn}

\begin{exa}
\label{E: Phi}
Assume $n = 6$.
Take $D \in \BPD([2,1,6,5,3,4])$
as depicted on the left. 
Then we perform the encoding algorithm and obtain $\Phi(D)$.

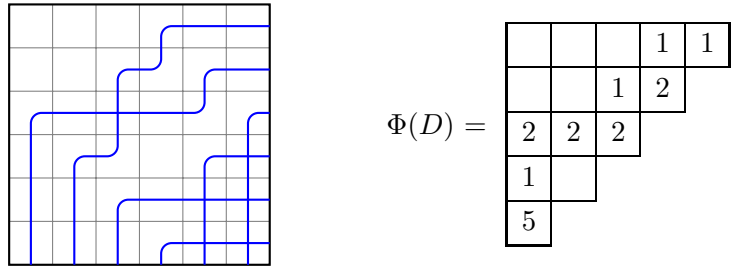
\[
\begin{tikzpicture}[x=1.5em,y=1.5em,thick,rounded corners, color = blue]
\draw[step=1,gray,thin] (0,0) grid (6,6);
\draw[color=black, thick, sharp corners] (0,0) rectangle (6,6);
\draw(0.5, 0)--(0.5,3.5)--(4.5,3.5)--(4.5,4.5)--(6,4.5);
\draw(1.5, 0)--(1.5,2.5)--(2.5,2.5)--(2.5,4.5)--(3.5,4.5)--(3.5,5.5)--(6,5.5);
\draw(2.5, 0)--(2.5,1.5)--(6,1.5);
\draw(3.5, 0)--(3.5,0.5)--(6,0.5);
\draw(4.5, 0)--(4.5,2.5)--(6,2.5);
\draw(5.5, 0)--(5.5,3.5)--(6,3.5);
\end{tikzpicture}
\quad\quad
\quad\quad
\raisebox{1.75cm}{$\Phi(D)$ = \raisebox{0.9cm}{
\begin{ytableau}
\: & \: & \: & 1 & 1\cr
\: & \: & 1 & 2\cr
2 & 2 & 2\cr
1 & \:\cr
5 \end{ytableau}}}
\]
\end{exa}

To precisely
describe the image of $\BPD(w)$ under $\Phi$,
we need the following definition.

\begin{defn}
The \definition{reading word} 
of a flagged tableau $T$, 
denoted as $\word(T)$,
is a sequence of pairs obtained as follows. 
Go through entries of $T$ from top to bottom, and right to left in each row.
When we see the number $i$ in  column $c$, 
we write the pair $(i, n+1-c)$.
\end{defn}

By the definition of flagged tableaux,
if we read $(i, j)$ in row $k$ of a flagged tableaux,
we must have $i \leq k < j$.

\begin{exa}
In Example~\ref{E: Phi}, 
$$\word(\Phi(D)) = 
(1,2)(1,3)(2,3)(1,4)(2,4)(2,5)(2,6)(1,6)(5,6).$$
\end{exa}

Let $T$ be a flagged tableau with reading word
$(a_1, b_1), \cdots, (a_d, b_d)$.
Then we say $T$ is \definition{associated} with the permutation
$w$ if 
$$
w \xrightarrow{t_{a_d, b_d}} \cdots \xrightarrow{t_{a_{2}, b_2}}
\xrightarrow{t_{a_{1}, b_1}}
w_0
$$
is a saturated Bruhat chain. Let $\FT(w)$
consist of all flagged tableaux
associated with $w$.
\begin{exa}
In Example~\ref{E: Phi},
$\Phi(D)$ is associated with 
$[2,1,6,5,3,4]$ because:
\begin{align*}
& [2,1,6,5,3,4] 
\xrightarrow{t_{5,6}} [2,1,6,5,4,3]
\xrightarrow{t_{1,6}} [3,1,6,5,4,2]
\xrightarrow{t_{2,6}} [3,2,6,5,4,1]
\xrightarrow{t_{2,5}} [3,4,6,5,2,1]\\
\xrightarrow{t_{2,4}} & [3,5,6,4,2,1]
\xrightarrow{t_{1,4}} [4,5,6,3,2,1]
\xrightarrow{t_{2,3}} [4,6,5,3,2,1]
\xrightarrow{t_{1,3}} [5,6,4,3,2,1]
\xrightarrow{t_{1,2}} [6,5,4,3,2,1]  
\end{align*}
is a saturated Bruhat chain
from $[2,1,6,5,3,4]$ to $w_0$.
Notice that $D \in \BPD([2,1,6,5,3,4])$.
\end{exa}

For a flagged tableau $T$,
define the \definition{weight} of $T$,
denoted as $\wt(T)$,
to be a sequence of $n-1$ numbers 
whose $i\textsuperscript{th}$
entry is the number of empty cells in row $i$. 
Then we have:

\begin{prop}
\label{P: BPD FT}
For $w \in S_n$, 
$\Phi$ is a weight-preserving bijection from 
$\BPD(w)$ to $\FT(w)$.
\end{prop}

The rest of this section
aims to prove Proposition~\ref{P: BPD FT}.
Recall from the previous section that 
we built a bijection $\chain_\BPD$
from $\BPD(w)$ to 
chains from $w$ to $w_0$ 
compatible with $(n-1, \cdots, 1)$.
The key step is to also encode these 
chains into flagged tableaux.
Let $C = (u_{n-1}, \cdots, u_1, u_0)$
be such a chain.
We construct a flagged tableau $\Psi(C)$
as follows.
Take $k \in [n-1]$.
Let $u_{k} \xrightarrow{t_{a_d, b_d}} \cdots \xrightarrow{t_{a_1, b_1}} u_{k-1}$
be the increasing $k$-chain from $u_k$ to $u_{k-1}$.
We know $a_i \leq k < b_i$ for each $i \in [d]$.
Then we place the number $a_i$
in row $k$ column $n + 1 - b_i$.
This cell exists since $n+1-b_i < n + 1 - k$
and row $k$ of $\Psi(C)$ has $n - k$ cells. 
The resulting filling is a flagged tableau
since $a_i \leq k$.

\begin{lemma}
\label{L: Psi}
The map $\Psi$ 
is a bijection from the set of 
chains from $w$ to $w_0$ 
compatible with $(n-1, \cdots, 1)$
to $\FT(w)$.
\end{lemma}
\begin{proof}
Let $C = (u_{n-1}, \cdots, u_1, u_0)$
be such a chain.   
Take $k \in [n-1]$.
Let $u_{k} \xrightarrow{t_{a_d, b_d}} \cdots \xrightarrow{t_{a_1, b_1}} u_{k-1}$
be the increasing $k$-chain from $u_k$ to $u_{k-1}$.
To show $\Psi(C) \in \FT(u_{n-1})$,
we just need to check the reading word
from row $k$ of $\Psi(C)$ 
is $(a_1, b_1), \cdots, (a_d, b_d)$.

Since the $k$-chain is increasing, 
$u_{k-1}(b_d) < \cdots < u_{k-1}(b_1)$.
By Corollary~\ref{C: chain in se},
we know $u_{k-1} \in S_n^{k \searrow}$,
so $b_1 < \cdots < b_d$.
Consequently, 
$n+1 - b_1 > \cdots > n+1-b_d$.
When we read row $k$ of $\Psi(C)$
from right to left, 
we see a cell in column $n+1-b_i$ 
filled by $a_i$,
where $i$ ranges from $1$ to $d$.
Thus, 
the reading word we get from
this row is $(a_1, b_1), \cdots, (a_d, b_d)$.

Finally, to show $\Psi$ is bijective,
we describe its inverse. 
Take $T \in \FT(w)$
and set $u_{0} = w_0$.
For $k \in [n-1]$,
suppose the reading word of row $k$ of $T$
is $(a_1, b_1), \cdots, (a_d, b_d)$.
We then define $u_k := u_{k-1} t_{a_1, b_1} \cdots t_{a_d, b_d}$.
Clearly, 
$u_{n-1} = w$.
Notice that $u_{k} \xrightarrow{t_{a_d, b_d}} \cdots \xrightarrow{t_{a_1, b_1}} u_{k-1}$
is an increasing $k$-chain if $u_{k-1} \in S_n^{k \searrow}$:
This is because $a_i \leq k < b_i$ for each $i$
and $b_1 < \cdots < b_d$ implies $u_{k-1}(b_1) > \cdots > u_{k-1}(b_d)$.
Finally, a simple induction using Corollary~\ref{C: chain in se}
would show $u_{k-1} \in S_n^{k \searrow}$
and the chain from $u_k$ to $u_{k-1}$ is indeed increasing.
Thus, $(u_{n-1}, \cdots, u_0)$ is compatible with $(n-1, \cdots, 1)$.
It is clear that this is the inverse of $\Psi$.
\end{proof}

The relation between the encoding $\Phi$ and the encoding $\Psi$
can be summarized as follows. 
\begin{lemma}
\label{L: Phi Psi}
Take $D \in \BPD(w)$.
We have $\Phi(D) = \Psi(\chain_\BPD(D))$.
\end{lemma}

\begin{proof}
Let $\chain_\BPD(D) = 
(u_{n-1}, \cdots, u_1, u_0)$.
Fix $k \in [n-1]$ and let
$(a_1, b_1), \cdots, (a_d, b_d)$
be the reading word coming from row 
$k$ of $\Phi(D)$.
We need to check the reading word
on row $k$ of $\Psi(\chain_\BPD(D))$
is also $(a_1, b_1), \cdots, (a_d, b_d)$.
Thus, it is enough to show
\begin{equation}
\label{eq: Phi proof}
u_k \xrightarrow{t_{a_d, b_d}}  \xrightarrow{t_{a_{d-1}, b_{d-1}}} \cdots \xrightarrow{t_{a_{2}, b_2}}
\xrightarrow{t_{a_{1}, b_1}}
u_{k-1}    
\end{equation}
is the increasing $k$-chain $u_k$ to $u_{k-1}$. 
First, we know there are $d$
$\rtile$ and $\htile$ in row $k$
of $D$.
Let $c_1 < \cdots < c_d$
be their column indices.
Then the pipe in $c_i$ is 
from row $a_i$.
Say $(k, c_i)$
is the $j_i\textsuperscript{th}$
tile in row $k$ that is 
$\rtile$, $\htile$ or $\btile$.
Then $b_i = n + 1 - j_i$.
By Corollary~\ref{c: One row to chain},
we know~\eqref{eq: chain from a row}
is an increasing $k$-chain.
\end{proof}

Now we are ready to prove Proposition~\ref{P: BPD FT}.
\begin{proof}[Proof of Proposition~\ref{P: BPD FT}]
Take $D \in \BPD(w)$.
By Lemma~\ref{L: Phi Psi} and Lemma~\ref{L: Psi}, 
$\Phi(D) \in \FT(w)$.
Then $\Phi(\cdot)$ is bijective
since $\Psi(\cdot)$ and $\chain_\BPD(\cdot)$ are.
Finally, it is clear
that the number of $\btile$ in row $k$ of $D$
agrees with the number of empty cells in
row $k$ of $\Phi(D)$.
\end{proof}

\subsection{Constructing  \texorpdfstring{$\fS^\BPD$}{}}
\label{SS: construction}
Define
$A:= \Q[x_1, \cdots, x_{n-1}] \otimes \E_n$.
It acts on $\Q[x_1, \cdots, x_{n-1}][S_n]$ from the right:
$(fw) \odot (g \otimes e) = (fg)(w \odot e)$
for any $f, g \in \Q[x_1, \cdots, x_{n-1}]$, $w \in S_n$
and $e \in \E_n$.
We may identify $\E_n$ and $\Q[x_1, \cdots, x_{n-1}]$
as subalgebras of $A$.
Now we define elements in $A$ that represent 
the possibilities of each cell and each row of
of a flagged tableau.

\begin{defn}
Take $i \in [n-1]$.
For $i < j$, define $B_{i, j} \in \E_n$
as 
$$B_{i,j} := 
d_{1, j} + \cdots + d_{i, j}.$$
Define $R_i(x_i) \in A$ as
$$R_i(x_i) := 
(x_i + B_{i, i+1})
(x_i + B_{i, i+2})
\cdots(x_i + B_{i, n}).$$
Finally, define 
$\fS^\BPD \in \Q[x_1, \cdots x_{n-1}][S_n]$ as 
$$\fS^\BPD := w_0 \odot 
(R_1(x_1) R_2(x_2) \cdots R_{n-1}(x_{n-1})).$$
\end{defn}

\begin{exa}
We compute $\fS^\BPD$ with $n = 3$.
We have 
$$
B_{1,2} = d_{1,2}, B_{1,3} = d_{1,3},\textrm{ and } 
B_{2,3} = d_{1,3} + d_{2,3}.
$$
Thus,
$$
R_1(x_1) = (x_1 + d_{1,2}) (x_1 + d_{1,3}) \textrm{ and }
R_2(x_2) = (x_2 +  d_{1,3} + d_{2,3}).
$$
We have 
\begin{align*}
\fS^\BPD & = [3,2,1] \odot 
(x_1 + d_{1,2})(x_1 + d_{1,3})(x_2 +  d_{1,3} + d_{2,3})\\
& = (x_1 [3,2,1] + [2,3,1]) \odot (x_1 + d_{1,3})(x_2 +  d_{1,3} + d_{2,3}) \\
& = (x_1^2 [3,2,1] + x_1[2,3,1] + [1,3,2]) \odot (x_2 +  d_{1,3} + d_{2,3}) \\
& = x_1^2x_2 [3,2,1] + x_1x_2[2,3,1] + (x_1 + x_2)[1,3,2]
+ x_1 [2,1,3] + x_1^2[3,1,2] + 1[1,2,3]
\end{align*}

\end{exa}

We show $\fS^\BPD$ is a generating function
of flagged tableaux, or equivalently all BPDs:

\begin{proof}[Proof of Proposition~\ref{P: generating function}]
If we expand $R_i(x_i)$, each term corresponds to
one way of filling row $i$ of a flagged tableau.
The expression $(x_i + B_{i, j})$ in $R_i(x_i)$
corresponds to ways of 
filling the cell at row $i$ 
and column $n+1-j$: 
$x_i$ means to leave the box empty
and $d_{p,j}$
means to fill it with $p$.
If we expand 
$R_1(x_1) \cdots R_{n-1}(x_{n-1})$,
for each term 
$x^\alpha d_{a_1, b_1} \cdots d_{a_k, b_k}$,
there is a flagged tableau 
$T$ with $\wt(T) = x^\alpha$
and $\word(T) = (a_1,b_1)\cdots (a_k, b_k)$.
Let $w = w_0 \odot d_{a_1, b_1} \cdots d_{a_k, b_k}$. 
If $w = 0$, we know $T$ is not associated with any permutation.
Otherwise, $T \in \FT(w)$.
Thus, we have the first equation.
The second equation follows from Proposition~\ref{P: BPD FT}. 
\end{proof}

\subsection{Proof of Theorem~\ref{T: FS on BPD}}
\label{SS: FS Proofs}
Our proof is similar to the arguments
of Fomin and Stanley. 
Consider a right action of $\N_n$ on $S_n$
with $w \odot u_i = w t_{i, i+1}$ if $w(i) < w(i+1)$
and $w \odot u_i = 0$ otherwise.
We may extend this action to $\Q[x_1, \cdots, x_{n-1}][S_n]$
by setting $f \odot u_i = f$ for all $f \in \Q[x_1, \cdots, x_{n-1}]$.
Similar to Fomin and Stanley's approach, 
Theorem~\ref{T: FS on BPD} reduces to: 
\begin{prop}
\label{P: FS on BPD}
For each $i \in [n-1]$,
\begin{align}
\label{EQ: FS on BPD}
\partial_i(\fS) = \fS \odot u_i.   
\end{align}
\end{prop}

Let us start with analyzing 
the left hand side of 
$(\ref{EQ: FS on BPD})$:
\begin{align}
\begin{split}
\label{EQ: FS LHS}
\partial_i(\fS) & = \partial_i(w \odot R_1(x_1) \cdots R_{n-1}(x_{n-1})) \\
& = w \odot \partial_i(R_1(x_1) \cdots R_{n-1}(x_{n-1})) \\
& = w \odot R_1(x_1) \cdots \partial_i(R_i(x_i)R_{i+1}(x_{i+1})) \cdots R_{n-1}(x_{n-1})     
\end{split}
\end{align}

To study the right hand side of 
$(\ref{EQ: FS on BPD})$,
we need the following lemma.

\begin{lemma}
\label{L: Swapping d and u}
Take $1 \leq i < i+1 < j \leq n$.
For any $a \in \Q[x_1, \cdots, x_n][S_n]$,
$$
a \odot (d_{i,j} u_{i}) = 
a \odot (u_{i} d_{i+1,j})
\quad\quad\quad\quad
a \odot (d_{i+1,j} u_{i}) = 
a \odot (u_{i} d_{i,j}).
$$
\end{lemma}
\begin{proof}
It is enough to assume $a = w$
for some $w \in S_n$.
The only case for $w \odot d_{i,j} u_i \neq 0$
is when $w(j) < w(i) < w(i+1)$
and $w(k) > w(i)$ or $w(k) < w(j)$
for any $i+1 < k < j$.
In this case, $w \odot d_{i,j} u_i = w t_{i,j} t_{i, i+1}$.
Correspondingly, the only case for $w \odot u_i d_{i+1,j} \neq 0$
is also when $w(j) < w(i) < w(i+1)$
and $w(k) > w(i)$ or $w(k) < w(j)$
for any $i+1 < k < j$.
In this case, $w \odot u_i d_{i+1,j} = w t_{i, i+1} t_{i+1, j} =  w t_{i,j} t_{i, i+1}$.
This proves the first equation. 
The second equation can be proved similarly. 
\end{proof}

\begin{lemma}
Take $1 \leq i < j \leq n$.
For any $a \in \Q[x_1, \cdots, x_n][S_n]$,
$$
a \odot R_j(x_j) u_{i}
= a \odot u_{i} R_j(x_j) 
$$
\end{lemma}
\begin{proof}
By Lemma~\ref{L: Swapping d and u}:
$$
a \odot (x_j + d_{1, k} + \cdots d_{j,k})u_{i}
= a \odot u_{i}(x_j + d_{1, k} + \cdots d_{j,k})
$$
for any $k > j$.
\end{proof}

Thus, right hand side of 
(\ref{EQ: FS on BPD})
becomes
\begin{align}
\label{EQ: FS RHS}
w_0 \odot
R_1(x_1) \cdots R_{i}(x_i)
\:u_{i, i+1}\: 
R_{i+1}(x_{i+1})
\cdots R_{n-1}(x_{n-1}).
\end{align}

Next, we want to describe 
the terms that 
appear in $w_0 \odot
R_1(x_1) \cdots R_{i-1}(x_{i - 1})$
with non-zero coefficients. 
Recall $S_n^{k \searrow}$ is the subset of $S_n$
consisting of all $w$ such that 
$w(k+1) > w(k + 2) > \cdots > w(n)$.
Let $\Q[S_n^{k \searrow}]$ be the subspace
of $\Q[S_n]$ with basis $S_n^{k \searrow}$.

\begin{lemma}
\label{L: Stay in S^k}
Take $k \in [n-1]$.
Then take $i,j \in [n-1]$ with $i \leq j$.
For any $a \in \Q[S_n^{k \searrow}]$, 
we have $a \odot d_{i, j} \in \Q[S_n^{k \searrow}]$.
\end{lemma}
\begin{proof}
Follows from Lemma~\ref{L: Stay in S^k (perm)}.
\end{proof}

\begin{cor}
\label{C: in smaller S}
For $k = 0, 1, \cdots, n-1$,
$w_0 \cdot R_1(x_1) \cdots R_k(x_k) \in 
\Q[x_1, \cdots, x_k][S_n^{k \searrow}]$.
\end{cor}
\begin{proof}
Prove by induction on $k$.
When $k = 0$, we have $w_0 \subseteq \Q[S_n^{0 \searrow}]$. 
Then the inductive step is finished by Lemma~\ref{L: Stay in S^k}.
\end{proof}

Combining Corollary~\ref{C: in smaller S}
with (\ref{EQ: FS LHS}) and (\ref{EQ: FS RHS}), 
Proposition~\ref{P: FS on BPD}
reduces to proving the following:
For any $i \in [n-1]$ and $a \in \Q[x_1, \cdots, x_{i-1}][S_n^{i-1 \searrow}]$,
\begin{align}
\label{EQ: FS Goal}
a \odot \partial_i(R_{i}(x_i) R_{i+1}(x_{i+1})) = a \odot R_{i}(x_i) u_i R_{i+1}(x_{i+1}).    
\end{align}

To prove this statement, we
introduce a powerful lemma
that turns $R_i$ into $R_{i+1}$ in some circumstances.
\begin{lemma}
\label{L: Break R_i}
For $i \in [n-2]$ and $a \in \Q[x_1, \cdots, x_n][S_n^{i\searrow}]$,
we have 
$$
a \odot R_i(y) = a \odot (y + \theta_{i+1}) R_{i+1}(y),
$$
where $y \in \{x_1, \cdots, x_n\}$.
Here, $\theta_{i+1}$
is the Dunkl element defined in Subsection~\ref{SS: FK}.
\end{lemma}
\begin{proof}
It is enough to assume $a = w$ for some $w \in S_n^{i\searrow}$.
We have
\begin{align}
\label{EQ: Break R_i}
\begin{split}
& w \odot (y + B_{i, i+1}) R_{i+1}(y)\\
= & w \odot (y + B_{i, i+1}) (y + B_{i+1, i+2}) (y + B_{i+1, i+3}) \cdots (y + B_{i+1}, n)\\
= & w \odot (y + B_{i, i+1}) (y + B_{i, i+2} + d_{i+1, i+2}) (y + B_{i+1, i+3}) \cdots (y + B_{i+1}, n)\\
= & w \odot (y + B_{i, i+1}) (y + B_{i, i+2}) (y + B_{i+1, i+3}) \cdots (y + B_{i+1}, n)\\
+ & w \odot (y + B_{i, i+1}) 
d_{i+1, i+2}(y + B_{i+1, i+3}) \cdots (y + B_{i+1}, n)\\
\end{split}
\end{align}

We have two summands in 
the last expression of (\ref{EQ: Break R_i}). 
Notice that for any $b \in \Q[x_1, \cdots, x_n][S_n^{i \searrow}]$, 
we have $b \odot B_{i+1, j} = b \odot B_{i,j}\in \Q[x_1, \cdots, x_n][S_n^{i \searrow}]$ if $j > i+2$,
so the first summand becomes
$$w \odot (y + B_{i, i+1}) (y + B_{i, i+2})(y + B_{i, i+3})\cdots (y + B_{i, n}) = w \odot R_i(y).$$
For the second summand, 
notice that:
\begin{align*}
& w \odot B_{i, i+1} d_{i+1, i+2} \\
= & w \odot \sum_{j = 1}^i d_{j, i+1} d_{i+1, i+2}\\
= & w \odot \sum_{j = 1}^i d_{i+1, i+2} d_{j, i+2}
+ w \odot \sum_{j = 1}^i d_{j, i+2} d_{j, i+1}\\
= & w \odot \sum_{j = 1}^i d_{i+1, i+2} d_{j, i+2}\\
= & w \odot d_{i+1, i+2} B_{i+1, i+ 2},
\end{align*}

where the second to last 
equality is given by 
$w \in \Q[x_1, \cdots, x_n][S_n^{i \searrow}]$.
All together, we have
\begin{align*}
w \odot (y + d_{1, i+1} + \cdots + d_{i, i+1}) R_{i+1}(y)
& = w \odot (y + B_{i, i+1}) R_{i+1}(y) \\
& =  w \odot R_i(y) + w \odot d_{i+1, i+2} R_{i+1}(y).
\end{align*}
After rearranging the terms, we have
$$
w \odot R_i(y) = 
w \odot (y + d_{1, i+1} + \cdots + d_{i, i+1} - d_{i+1, i+2})R_{i+1}(y).
$$

To see the right hand side
is $w \odot (y - \theta_{i+1}) R_{i+1}(y)$,
just notice that
$w \odot d_{i+1, j} = 0$
for any $j > i+2$.
\end{proof}
\begin{rem}
To make Lemma~\ref{L: Break R_i} still work for $i = n-1$,
we simply 
define $R_n(y)$ as $1$.
Then it is easy to see
$R_{n-1}(y) = (1 - \theta_{n}) R_n(y)$.
\end{rem}

Then we have the following 
commutativity of Dunkl
elements and $R_i(y)$.

\begin{lemma}
\label{L: theta R}
For $i, j \in [n]$ with $i \leq j$ and $a \in \Q[x_1, \cdots, x_n][S_n^{j\searrow}]$,
$a \odot \theta_i R_j(y) = a \odot R_j(y) \theta_i$,
where $y \in \{x_1, \cdots, x_n\}$.
\end{lemma}
\begin{proof}
Prove by induction on $j$.
When $j = n$, the lemma
is immediate.
Now if $j < n$,
we have $a \odot \theta_i 
\subset \Q[x_1, \cdots, x_n][S_n^{j \searrow}]$ by $i \leq j$.
Lemma~\ref{L: Break R_i} says 
$$
a \odot \theta_i R_j(y) = a \odot \theta_i (y + \theta_{j+1}) R_{j+1}(y).$$
Since the Dunkl elements commute, 
we have 
$$
a \odot \theta_i R_j(y) 
= a \odot (y + \theta_{j+1}) \theta_i R_{j+1}(y)
= a \odot (y + \theta_{j+1}) R_{j+1}(y)\theta_i= a \odot R_{j}(y)\theta_i$$
where the second equation is by our inductive hypothesis 
and the last equation is by Lemma~\ref{L: Break R_i}.
\end{proof}

The commutativity of Dunkl elements and $R_i(y)$
allows us to derive the 
following symmetry.
\begin{lemma}
\label{L: R symmetry}
Take $i\in [n]$ and $a \in \Q[x_1, \cdots, x_n][S_n^{i\searrow}]$.
Consider distinct $y, z \in \{x_1, \cdots, x_n\}$.
Assume $a$ is symmetric
in $y, z$,
then so is
$a \odot R_{i}(y) R_{i}(z)$.
\end{lemma}
\begin{proof}
Prove by induction on $i$.
When $i = n$,
the statement is immediate.
Assume $i < n$.
By Lemma~\ref{L: Break R_i}
and Lemma~\ref{L: theta R},
we have
\begin{align*}
a \odot R_{i}(y) R_{i}(z)
= & a \odot R_{i}(y) (z - \theta_{i+1}) R_{i+1}(z)\\
= & a \odot R_{i}(y) R_{i+1}(z)(z - \theta_{i+1})\\
= & a \odot (y - \theta_{i+1}) R_{i+1}(y) R_{i+1}(z)(z - \theta_{i+1})\\
= & a \odot  R_{i+1}(y) R_{i+1}(z)(y - \theta_{i+1})(z - \theta_{i+1})\\
\end{align*}

Our inductive hypothesis says $a \odot  R_{i+1}(y) R_{i+1}(z)$
is symmetric in $y$ and $z$,
so the inductive step is finished. 
\end{proof}

Finally, we are ready to prove Proposition~\ref{P: FS on BPD}.
\begin{proof}[Proof of Proposition~\ref{P: FS on BPD}]
We just need to check (\ref{EQ: FS Goal}) holds 
for any $a \in \Q[x_1, \cdots, x_{i-1}][S_n^{i-1 \searrow}]$ and any
$i \in [n-1]$.
By Lemma~\ref{L: Break R_i} and Lemma~\ref{L: theta R},
the left hand side of (\ref{EQ: FS Goal}) is:
\begin{align*}
& \partial_i(a \odot R_{i}(x_i) R_{i+1}(x_{i+1}))\\
= & \partial_i(a \odot (x_i - \theta_{i+1}) R_{i+1}(x_i) R_{i+1}(x_{i+1})) \\
= & \partial_i(a \odot R_{i+1}(x_i) R_{i+1}(x_{i+1})(x_i - \theta_{i+1}) ) \\
= & (a \odot R_{i+1}(x_i) R_{i+1}(x_{i+1})) \odot \partial_i(x_i - \theta_{i+1} ) \\
= & a \odot R_{i+1}(x_i) R_{i+1}(x_{i+1}),
\end{align*}
where the second to last equation is by Lemma~\ref{L: R symmetry}.
For the right hand side of (\ref{EQ: FS Goal}),
by Lemma~\ref{L: Break R_i} and Lemma~\ref{L: theta R}:
\begin{align*}
& a \odot R_{i}(x_i) u_{i} R_{i+1}(x_{i+1})\\
= & a \odot (x_i - \theta_{i+1}) R_{i+1}(x_i) u_{i} R_{i+1}(x_{i+1})) \\
= & a \odot R_{i+1}(x_i) (x_i - \theta_{i+1}) u_{i} R_{i+1}(x_{i+1})) \\
= & a \odot R_{i+1}(x_i)  R_{i+1}(x_{i+1})) \\
\end{align*}

The last equality needs some explanations. 
For any permutation $w$
that has non-zero coefficient in
$a \odot R_{i+1}(x_i)$,
we know $w(i) > w(i+1)$.
Then $w u_i = 0$ and $w d_{j, i+1} u_i = w d_{i+1, j} u_i = 0$
if $j \neq i$.
Notice that $w d_{i, i+1}, u_i = w$,
so $w (x_i - \theta_{i+1}) u_i = w$.
\end{proof}

\subsection{Commutative expression of  \texorpdfstring{$\fS^\BPD$}{}}
\label{SS: Commuative}
We end this subsection by providing an alternative way to write $\fS^\BPD$ as a commutative product using Dunkl elements.

The lemmas we derived along the
way leads to an alternative formula of $\fS$:

\begin{prop}
\label{P: Commutative FS}
We have:
\begin{align}
\label{EQ: ultimate}
\fS = w_0 \odot \prod_{1 \leq i < j \leq n} (x_i - \theta_j).
\end{align}
\end{prop}
Notice that terms 
we multiply on the right
hand side commute with each other, so may use the $\prod$
notation. 
\begin{proof}
We claim for any $k = 0, \cdots, n-2$ and  
$a \in \Q[x_1, \cdots, x_n][S_n^{k \searrow}]$,
$$
a \odot R_{k+1}(x_{k+1}) \cdots R_{n-1}(x_{n-1})
= a \odot \prod_{k+1 \leq i < j \leq n} (x_i - \theta_j).
$$
Then the theorem can be obtained
by setting $k = 0$.
We prove our claim by induction
on $k$.
When $k = n-2$,
just observe that $R_{n-1}(x_{n-1}) = (x_{n-1} - \theta_n)$.
Now assume $k < n-2$.
We just keep using 
Lemma~\ref{L: Break R_i} and Lemma~\ref{L: theta R}:
\begin{align*}
& a \odot R_{k+1}(x_{k+1}) R_{k+2}(x_{k+2})\cdots R_{n-1}(x_{n-1}) \\
= & a \odot R_{k+1}(x_{k+1}) 
\prod_{k+2 \leq i < j \leq n} (x_i - \theta_j)\\
= & a \odot R_{k+2}(x_{k+1})(x_{k+1} - \theta_{k+2})\prod_{k+2 \leq i < j \leq n} (x_i - \theta_j)\\
= & a \odot R_{k+3}(x_{k+1})(x_{k+1} - \theta_{k+3})(x_{k+1} - \theta_{k+2})\prod_{k+2 \leq i < j \leq n} (x_i - \theta_j)\\
= & \cdots \\
= & \prod_{k+1 \leq i < j \leq n} (x_i - \theta_j) \qedhere
\end{align*}
\end{proof}

\begin{exa}
In this example, 
we compute $\fS$ using (\ref{EQ: ultimate})
for $n = 3$.
\begin{align*}
\fS = & [3,2,1] \odot (x_1 - \theta_2) (x_1 - \theta_3) (x_2 - \theta_3) \\
= & (x_1[3,2,1] + [2,3,1] - [3,1,2]) \odot (x_1 - \theta_3) (x_2 - \theta_3)\\
= & (x_1^2[3,2,1] + x_1[2,3,1] - x_1[3,1,2] + x_1[3,1,2] + [1,3,2] + [2,1,3] - [2,1,3]) \odot (x_2 - \theta_3)\\
= & (x_1^2[3,2,1] + x_1[2,3,1] + [1,3,2] ) \odot (x_2 - \theta_3)\\
= & x_1^2x_2[3,2,1] + x_1x_2[2,3,1] + x_2[1,3,2] + x_1^2[3,1,2] + x_1[1,3,2] + x_1[2,1,3] + [1,2,3]\\
= & x_1^2x_2[3,2,1] + x_1x_2[2,3,1] + (x_1 + x_2)[1,3,2] + x_1^2[3,1,2] + x_1[2,1,3] + [1,2,3]
\end{align*}
\end{exa}

\begin{rem}
Sergey Fomin kindly informed the author the term
$w_0 \odot \prod_{1 \leq i < j \leq n} (x_i - \theta_j)$ in~(\ref{EQ: ultimate})
seems related to the following variation of Cauchy identity of Schubert polynomials:
\begin{align}
\label{EQ: Cauchy}
\prod_{1 \leq i < j \leq n} (x_i - y_j) = \sum_{w \in S_n} \fS_w(x_1, \cdots, x_{n-1})\fS_{ww_0}(-y_n, \cdots, -y_2).
\end{align}
Indeed, there is a relation. 
Following Lemma~\ref{L: Dunkl Monk},
if we replace $y_j$ in~\eqref{EQ: Cauchy}
by $\theta_j$ and apply it on $w_0$,
the equation becomes $w_0 \odot \prod_{1 \leq i < j \leq n} (x_i - \theta_j) = \sum_{w \in S_n}\fS_w w$.
Notice that the two sides are the expressions in 
Proposition~\ref{P: Commutative FS}
and Theorem~\ref{T: FS on BPD}.
\end{rem}

\section{Bijection between pipedreams and bumpless pipedreams}
\label{S: Growth}

Our second application is a weight preserving
bijection between $\PD(w)$ and $\BPD(w)$. 
By Theorem~\ref{T: PD to chain} and Theorem~\ref{T: BPD to chain},
we just need a weight reversing bijection
between chains from $w$ to $w_0$ compatible with $(1,\cdots, n-1)$ and those compatible with $(n-1,\cdots, 1)$. 
This task can be done by Lenart's growth 
diagram~\cite{L_growth}.

\subsection{Lenart's growth diagram and a partial solution of
Problem~\ref{prob: chain}}

Given a saturated Bruhat chain 
$(w_1,\cdots, w_d)$,
we say its length is $d-1$.
We describe Lenart's growth diagram as the following algorithm.

\begin{defn}[\cite{L_growth}]
\label{D: growth move}
Given $k_1, k_2 \in [n-1]$
and chains $C_1, C_2$,
where $C_1$ (resp. $C_2$)
is a saturated $k_1$-chain from $u$ to $v$
(resp. saturated $k_2$-chain from $v$ to $w$).
\definition{Lenart's growth diagram}
is an algorithm that takes $k_1, k_2, C_1, C_2$
as inputs. 
It outputs a saturated $k_2$-chain from $u$ to some $v'$ and a 
saturated $k_1$-chain from $v'$ to $w$.
Moreover, the $k_1$-chain 
(resp. $k_2$-chain) in the output
has the same length as $C_1$ (resp. $C_2$). 

Assume $C_1 = (u_1, \cdots, u_{d_1})$ and 
$C_2 = (w_1, \cdots, w_{d_2})$ where $u_{d_1} = w_1$.
We first draw:
$$
u_1 \xrightarrow{k_1} u_2 \xrightarrow{k_1}
\cdots \xrightarrow{k_1} u_{d_1 - 1} \xrightarrow{k_1} w_1
\xrightarrow{k_2} w_2 \xrightarrow{k_2}
\cdots \xrightarrow{k_2} w_{d_2}.
$$
We start from this labeled chain and apply 
\definition{Lenart's local move}.
This move is applied to 
a segment of the chain
that looks like $a \xrightarrow{k_1} b \xrightarrow{k_2} c$.
We must have $a \lessdot_{k_1} b \lessdot_{k_2} c$.
There exists a unique $b' \in S_n$
such that $b' \neq b$
and $a \lessdot b' \lessdot c$.
If $a \lessdot_{k_2} b' \lessdot_{k_1} c$,
the local move will replace this part of the chain by
$a \xrightarrow{k_2} b' \xrightarrow{k_1} c$.
Otherwise, 
we must have $a \lessdot_{k_2} b' \lessdot_{k_1} c$.
The local move replace this part
by $a \xrightarrow{k_2} b \xrightarrow{k_1} c$.
We keep applying Lenart's local move until 
the labeled chain looks like:
$$
u_1' \xrightarrow{k_2} u_2' \xrightarrow{k_2}
\cdots \xrightarrow{k_2} u_{d_2 - 1}' \xrightarrow{k_2} w_1'
\xrightarrow{k_1} w_2' \xrightarrow{k_1}
\cdots \xrightarrow{k_1} w_{d_1}'.
$$
Then we output
the $k_2$-chain $(u_1', \cdots, u_{d_2-1}', w_1')$ and the $k_1$-chain
$(w_1', \cdots, w_{d_1}')$.
The ordering of the local moves clearly will not affect the output. 
\end{defn}

\begin{exa}
\label{E: growth}
Say the inputs are: $k_1 = 2$, $k_2 = 3$, 
$C_1 = ([2,1,4,3], [2,4,1,3], [3,4,1,2])$, 
and $C_2 = ([3,4,1,2], [3,4,2,1])$.
We start from the following labeled chain
and apply local moves: 
$$
[2,1,4,3] \xrightarrow{2} [2,4,1,3] \xrightarrow{2} [3,4,1,2]
\xrightarrow{3} [3,4,2,1].
$$
$$
[2,1,4,3] \xrightarrow{2} [2,4,1,3] \xrightarrow{3} [2,4,3,1]
\xrightarrow{2} [3,4,2,1],
$$
$$
[2,1,4,3] \xrightarrow{3} [2,3,4,1] \xrightarrow{2} [2,4,3,1]
\xrightarrow{2} [3,4,2,1].
$$
Therefore, the outputs are $([2,1,4,3], [2,3,4,1])$ and $([2,3,4,1], [2,4,3,1], [3,4,2,1])$. 
\end{exa}

The Lenart's local move is reversible. 
Assume it changes $a \xrightarrow{k_1} b \xrightarrow{k_2} c$
into $a \xrightarrow{k_1} b' \xrightarrow{k_2} c$.
Then applying it on $a \xrightarrow{k_1} b' \xrightarrow{k_2} c$
would yield $a \xrightarrow{k_1} b \xrightarrow{k_2} c$.
Therefore, Lenart concluded the following.
\begin{thm}\cite{L_growth}
\label{T: Lenart Bijective}
Suppose you input $k_1, k_2, C_1, C_2$ to Lenart's growth diagram
and obtain $C_2', C_1'$.
Then if you input $k_2, k_1, C_2', C_1'$ to Lenart's growth diagram,
you would obtain $C_1, C_2$ back. 
\end{thm}

As a consequence, for $u, w \in S_n$ and $k_1, k_2 \in [n-1]$,
Lenart's growth diagram gives a bijection 
between the following two sets: 
\begin{itemize}
\item The set of $(C_1, C_2)$
where $C_1$ is a saturated $k_1$-chain from $u$
to some $v$ and $C_2$ is a saturated $k_2$-chain from $v$ to $w$.
\item The set of $(C_2, C_1)$
where $C_2$ is a saturated $k_2$-chain from $u$
to some $v$ and $C_1$ is a saturated $k_1$-chain from $v$ to $w$.
\end{itemize}

Readers might expect that Lenart's growth 
diagram solves Problem~\ref{prob: chain}.
More explicitly, 
the following map seems 
like a solution.
\begin{defn}
Take a chain $(u,v,w)$ that is compatible with
$(k_1, k_2)$.
Let $C_1$ (resp. $C_2$)
be the increasing $k_1$-chain (resp. $k_2$-chain)
from $u$ to $v$ (resp. $v$ to $w$).
Input $k_1, k_2, C_1, C_2$
to Lenart's growth diagram,
obtaining a $k_2$-chain from $u$
to $v'$ and a $k_1$-chain from $v'$
to $w$.
Then define $\growth_{k_1, k_2}(u,v,w)$
as $(u, v', w)$.
\end{defn}

Unfortunately, this is not the solution. 
Notice that the input chains are both increasing chains
in Example~\ref{E: growth}.
However, the output $2$-chain is not increasing. 
There is no increasing $2$-chain from 
$[2,3,4,1]$ to $[3,4,2,1]$.
Therefore, the map $\growth_{2,3}$
would send $([2,1,4,3], [2,4,1,3],[3,4,1,2])$
to $([2,1,4,3], [2,3,4,1], [3,4,1,2])$
which is not compatible with $(3,2)$.
Nevertheless, in Corollary~\ref{C: growth 2},
we are going to show $\growth_{k_1, k_2}$
solves problem~\ref{prob: chain} in a special case.
The key step is to establish the following. 

\begin{lemma}
\label{L: Growth preserves inc}
Take $k \in [n-1]$ and $w \in S_n^{k \searrow}$.
Let $C_1$ be an increasing $k$ chain from $u$ to $v$
and $C_2$ be an increasing $1$ chain from $v$ to $w$.
Input $(k, 1, C_1, C_2)$ to Lenart's growth diagram,
obtaining $1$-chain $C_2'$ from $u$ to some $v'$
and $k$-chain $C_1'$ from $v'$ to $w$.
Then we know $C_2'$ is an increasing $1$-chain
and $C_1'$ is an increasing $k$-chain.
\end{lemma}

The rest of this section aims to prove this lemma. 
We first want to understand how Lenart's local move
behave in this case. 

\begin{rem}
\label{R: Local move}
Take $k \in [n-1]$ and $w \in S_n^{k \searrow}$.
We would like to describe how to apply 
Lenart's local move to a segment that looks like
$$
u \xrightarrow{k} v \xrightarrow{1} w.
$$
We may write $v = w t_{1,a}$
and $u = v t_{b,c}$
with $1 < a$ and $b \leq k < c$.
We explicitly state what Lenart's local move does in this case:
\begin{itemize}
\item Case 1: $1, a, b, c$ are distinct.
This segment becomes
$$
u \xrightarrow{1} w t_{b,c} \xrightarrow{k} w.
$$  
\item Case 2: $a = b$ and $w(b) > w(c)$.
We must have $w(1) > w(b) > w(c)$.
Notice that $w t_{b,c} \lessdot_{k} w$
and $u = w t_{b,c} t_{1,c} \lessdot_1 w t_{b,c}$.
This segment becomes
$$
u \xrightarrow{1} w t_{b,c} \xrightarrow{k} w.
$$  

\item Case 3: $a = b$ and $w(b) < w(c)$.
We must have $w(1) > w(c) > w(b)$.
Notice that $w t_{1,c} \lessdot_{k} w$
and $u = w t_{1,c} t_{1,b} \lessdot_1 w t_{1,c}$.
This segment becomes
$$
u \xrightarrow{1} w t_{1,c} \xrightarrow{k} w.
$$  

\item Case 4: $a = c$.
We must have $w(b) > w(1) > w(c)$.
Notice that $w t_{b,c} \lessdot_{k} w$
and $u = w t_{b,c} t_{1,b} \lessdot_1 w t_{b,c}$.
This segment becomes
$$
u \xrightarrow{1} w t_{b,c} \xrightarrow{k} w.
$$  

\item Case 5: $b  = 1$.
We must have $w(1) > w(a) > w(c)$.
This implies $a < c$.
If not, we have $k < c < a$ and $w(a) > w(c)$,
contradicting to  $w \in S_n^{k \searrow}$.

Now, notice that $w t_{a,c} \lessdot w$
and $u = w t_{a,c} t_{1,a} \lessdot_1 w t_{a,c}$.
However, we might not have $w t_{a,c} \lessdot_k w$.
Therefore, we need two sub-cases in this case.
The Case 5A corresponds to $a \leq  k$,
so $w t_{a,c} \lessdot_k w$.
This segment becomes
$$
u \xrightarrow{1} w t_{a,c} \xrightarrow{k} w.
$$  

Otherwise, the Case 5B corresponds to $k < a$.
This segment remains unchanged: 
$$
u \xrightarrow{1} w t_{1,a} \xrightarrow{k} w.
$$  
\end{itemize}
\end{rem}

\begin{lemma}
\label{L: 5B}
In Case 5B of Remark~\ref{R: Local move}, we can deduce 
$w(a) = w(c) + 1$ and $a = c - 1$.
\end{lemma}
\begin{proof}
When we swap the first number and $c\textsuperscript{th}$
number in $w t_{1,a}$,
we are swapping the values $w(a)$ 
and $w(c)$.
By $w t_{1,a}t_{1,c} = u \lessdot w t_{1,a}$,
if there is a value strictly between $w(a)$ and $w(c)$,
it cannot live on the left of $w(c)$ in $w t_{1,a}$.
On the other hand, recall $w \in S_n^{k \searrow}$.
By Lemma~\ref{L: Stay in S^k (perm)},
$w t_{1,a} \in S_n^{k \searrow}$.
A value larger than $w(c)$
cannot appear after $w(c)$ in $w t_{1,a}$ since $k < a < c$.
We conclude there is no value between $w(a)$ and $w(c)$,
so $w(a) = w(c) + 1$.
Then $a = c-1$ follows from $w \in S_n^{k \searrow}$.
\end{proof}

We are ready to prove Lemma~\ref{L: Growth preserves inc}
\begin{proof}[Proof of Lemma~\ref{L: Growth preserves inc}]
Let $C_1$ be $(u_1, \cdots, u_d)$ where $u = u_1$.
We first study the special case where chain $C_2$ has length $1$.
We may assume $C_2 = (u_d, w_d)$ where $w_d = w \in S_n^{k \searrow}$.
The algorithm starts with the picture 
$$
u_1 \xrightarrow{k} u_2 \xrightarrow{k} \cdots \xrightarrow{k} u_{d-1} \xrightarrow{k} u_{d} \xrightarrow{1} w_d.
$$
Then the algorithm applies Lenart's local move.
After applying $d-i$ moves, the picture looks like 
$$
u_1 \xrightarrow{k} u_2  \xrightarrow{k} \cdots \xrightarrow{k} u_{i}
\xrightarrow{1} w_{i} \xrightarrow{k} w_{i+1} \xrightarrow{k} \cdots \xrightarrow{k} w_{d}.
$$
Eventually, it becomes
$$
u_1 \xrightarrow{1} w_{1} \xrightarrow{k} w_{2} \xrightarrow{k} \cdots \xrightarrow{k} w_{d},
$$
so the output will be $(u_1, w_1)$
and $(w_1, \cdots, w_d)$.
Our goal is to show $(w_1, \cdots, w_d)$
is an increasing $k$-chain. 

Suppose $u_i = u_{i+1} t_{b_i, c_i}$
and $w_i = w_{i+1} t_{b_i', c_i'}$.
Since $w_d \in S_n^{k \searrow}$,
Lemma~\ref{L: Stay in S^k (perm)} implies
$u_1, \cdots, u_d$ and $w_1, \cdots, w_d$
are all in $S_n^{k \searrow}$.
Our goal is equivalent to show 
$c_1' > c_2' > \cdots > c_{d-1}'$.
Since $C_1$ is an increasing $k$-chain,
we know $c_1 > c_2 > \cdots > c_{d-1}$.

Now we consider one move of the algorithm.
It change a segment that looks like
$u_i \xrightarrow{k} u_{i+1} \xrightarrow{1} w_{i+1}$
to $u_i \xrightarrow{1} w_{i} \xrightarrow{k} w_{i+1}$.
Since $w_{i+1} \in S_n^{k \searrow}$,
this move's behavior is captured by Remark~\ref{R: Local move}.
We conclude that $c_i' = c_{i}$ except in Case 5B.
In this case, by Lemma~\ref{L: 5B}, $c_i' = c_i - 1$.
We need to make sure $c_i' > c_{i+1}'$ if $i < d - 1$.
Suppose not, we must have $c_{i+1}' = c_{i}' = c_i - 1$.
By Lemma~\ref{L: 5B},
if $w_{i+1}(c_i) = m$,
then $w_{i+1}(c_i-1) = m+1$.
By $c_{i+1}' = c_i - 1$,
$w_{i+2}(c_i-1) 
< w_{i+1}(c_i-1) = m+1$, 
so it is less than $m$
since $w_{i+2}(c_i) = m$.
This contradicts to $w_{i+2} \in S_n^{k \searrow}$.

Now we study the general case of 
Lemma~\ref{L: Growth preserves inc}:
$C_2$ might not have length $1$.
The algorithm starts from a picture
that looks like:
$$
u_1 \xrightarrow{k} u_2 \xrightarrow{k} \cdots \xrightarrow{k} u_{d-1} \xrightarrow{k} u_{d} \xrightarrow{1} w_d
\xrightarrow{1} w_{d+1}
\xrightarrow{1} \cdots 
\xrightarrow{1} w_D.
$$
When we apply Lenart's local moves
and swap a $1$ arrow all the way left, the consecutive $k$-chain
remains an increasing $k$-chain
by the special case above. 
Thus, the output $k$-chain of 
the algorithm is an increasing $k$-chain. 
Then the lemma is finished 
by noticing that any $1$-chain
is increasing.
\end{proof}

\begin{cor}
\label{C: growth 1}
Take $k \in [n-1]$,
$u \in S_n$
and $w \in S_n^{k \searrow}$.
The map $\growth_{k, 1}$
is a weight-reversing bijection from 
chains $(u,v,w)$ compatible with $(k,1)$
to chains $(u, v', w)$
compatible with $(1,k)$.
Its inverse is $\growth_{(1,k)}$.
\end{cor}
\begin{proof}
Take $(u,v,w)$ compatible with $(k,1)$.
Let $C_1$ be the increasing $k$-chain
from $u$ to $v$ and $C_2$
be the increasing $1$-chain from $v$ to $w$.
By Lemma~\ref{L: Growth preserves inc},
applying Lenart's growth diagram on
$k,1, C_1, C_2$ would yield
$C_2', C_1'$ which are increasing $1$-chain
and increasing $k$-chain respectively. 
Say $C_2'$ is from $u$ to $v'$,
so $\growth_{k,1}(u,v,w) = (u, v', w)$,
which is compatible with $(1,k)$.
Moreover, by Theorem~\ref{T: Lenart Bijective},
$\growth_{1,k}(u,v',w) = (u, v, w)$
and $\growth_{1,k}$ is weight-reversing. 

Now we have shown $\growth_{k,1}$
is an injective map from 
$(u,v,w)$ compatible with $(k,1)$
to $(u, v', w)$ compatible with $(1,k)$.
The proof is finished by Proposition~\ref{P: Swap chains},
which says there are the same number of
such $(u,v,w)$ and $(u, v', w)$.
\end{proof}

We end this subsection by 
slightly extending Corollary~\ref{C: growth 1}.
Our extension needs the following
notion.
\begin{defn}
\label{D: se k}
Let $S_n^{\searrow k}$ be the set of 
permutation $w \in S_n$ 
such that $w(i) = n + 1 - i$
for $i \in [k-1]$.
\end{defn}

\begin{cor}
\label{C: growth 2}
Take $k_1, k_2 \in [n-1]$
with $k_1 < k_2$.
Take $w \in S_n^{k_2 \searrow}$
and $u \in S_n^{\searrow k_1}$.
The map $\growth_{k_2, k_1}$
is a weight-reversing bijection from 
chains $(u,v,w)$ compatible with $(k_2,k_1)$
to chains $(u, v', w)$
compatible with $(k_1,k_2)$.
Its inverse is $\growth_{(k_1,k_2)}$.
\end{cor}
\begin{proof}
The condition guarantees the first $k_1 - 1$
numbers of $u$ are $n, \cdots, n - k_1 + 2$.
They will stay unswapped in an increasing $k_2$-chain,
so $v \in S_n^{\searrow k_1}$.
Similarly, the first 
$k_1 - 1$ numbers of $v$ will stay unswapped 
in an increasing $k_1$-chain, 
so $w \in S_n^{\searrow k_1}$.
By ignoring the first $k_1 - 1$ entries,  
we may view $u, v, w$ as $\tilde{u}, \tilde{v}, \tilde{w}$
in $S_{n - k_1 + 1}$.
Then $(\tilde{u}, \tilde{v}, \tilde{w})$
is compatible with $(k_2 - k_1 + 1, 1)$
where $\tilde{w} \in S_{n - k_1 + 1}^{k_2 - k_1 + 1 \searrow}$.
The proof is finished by Corollary~\ref{C: growth 1}
\end{proof} 

\subsection{Bijection between PDs and BPDs via Lenart's growth diagram.}

Now we use the map $\growth_{k_1,k_2}$ to describe a bijection 
from $\BPD(w)$ to $\PD(w)$.
We start with a simple lemma that
will be used to show our bijection is well-defined.
\begin{lemma}
\label{L: stay in se k}
Suppose there is a chain from $w$ 
to $w_0$ compatible with $(k_1, \cdots, k_d)$
where $k_1, \cdots, k_d$ are distinct numbers
from $[n-1]$.
Let $k$ be the smallest among them. 
Then $w \in S_n^{\searrow k}$.
\end{lemma}
\begin{proof}
By Proposition~\ref{P: Swap chains},
it is enough to assume $k_1 < \cdots < k_d$.
Then we prove by induction on $d$.
The lemma is immediate when $d = 1$.
Now suppose $d > 1$.
There is an increasing $k_1$-chain from 
$w$ to $w'$ and a chain from $w'$ 
to $w_0$ compatible with $(k_2, \cdots, k_d)$.
By our inductive hypothesis, since $k_2 > k_1$,
$w'(i) = n + 1 - i$ for $i \in [k]$.
The increasing $k$-chain from $w$ to $w'$
cannot move the first $k-1$ numbers of $w$,
so $w(i) = n+1 - i$ for $i \in [k-1]$.
\end{proof}

\begin{defn}
\label{D: BPD PD Bijection}
Let $C$ be a chain from $w$ to $w_0$ 
that is compatible with $(n-1, \cdots, 1)$.
We describe a map $\flip$ that sends it to a chain 
from $w$ to $w_0$ compatible with $(1, \cdots, n-1)$.

We will recursively compute chains $C_\gamma$,
where $\gamma$ range over certain permutations in $S_{n-1}$.
Moreover, each $C_\gamma$
is a chain from $w$ to $w_0$ compatible with $\gamma$.
For the base case, we set $C_{[n-1, \cdots, 1]} = C$.
Suppose we have $C_\gamma = (w_1, \cdots, w_n)$
where $\gamma$ is not the identity. 
Let $i$ be the smallest such that $\gamma(i) \neq i$.
Say $\gamma(a) = i$,
so $(w_{a-1}, w_a, w_{a+1})$ is compatible
with $(\gamma(a-1), i)$.
Let $\growth_{(\gamma(a-1), i)}(w_{a-1}, w_a, w_{a+1}) 
= (w_{a-1}, w_a', w_{a+1})$.
Then define
$$C_{\gamma t_{a-1, a}} = (w_1, \cdots,w_{a-1}, w_a', w_{a+1},\cdots, w_n).$$
Eventually, we will obtain $C_{[1,2, \cdots, n-1]}$,
which we define as $\flip(C)$.
\end{defn}

\begin{exa}
Suppose $n = 4$ and consider the chain 
$$C = ([2,1,4,3],[2,3,4,1],[2,4,3,1],[4,3,2,1])$$
which is
compatible with $(3,2,1)$.
We show how to compute $\flip(C)$, 
which is $C_{[1,2,3]}$.
The chain $C_{[3,2,1]}$ is set as $C$.
By $\growth_{2,1}([2,3,4,1],[2,4,3,1],[4,3,2,1]) = ([2,3,4,1],[4,2,3,1],[4,3,2,1])$,
we have 
$$
C_{[3,1,2]} = ([2,1,4,3],[2,3,4,1],[4,2,3,1],[4,3,2,1]).
$$
Since $\growth_{3,1}([2,1,4,3],[2,3,4,1],[4,2,3,1]) = ([2,1,4,3],[4,1,3,2],[4,2,3,1])$,
we have
$$
C_{[1,3,2]} = ([2,1,4,3],[4,1,3,2],[4,2,3,1],[4,3,2,1]).
$$
Finally, since $\growth_{3,2}([4,1,3,2],[4,2,3,1],[4,3,2,1]) = ([2,1,4,3],[4,3,1,2],[4,2,3,1])$, we have
$$
\flip(C) = C_{[1,2,3]} = ([2,1,4,3],[4,1,3,2],[4,3,1,2],[4,3,2,1]).
$$
\end{exa}

\begin{rem}
\label{R: gamma}
The set of $\gamma$ that shows up 
in the definition of $\flip$
can be described as follows. 
They are permutations obtained by inserting
a number $i \in [n-1]$ to 
$(n-1, \cdots, i+1)$,
and then prepend $(1,\cdots, i-1)$
in front. 
For example, 
when $n = 4$,
we list the $\gamma$ that appear in order:
$$
[4,3,2,1], [4,3,1,2],
[4,1,3,2], [1,4,3,2],
[1,4,2,3], [1,2,4,3],
[1,2,3,4].
$$
Let $\gamma$ be such a permutation in $S_{n-1}$
that is not $[1,\cdots, n-1]$.
Let $i$ be the smallest such that $\gamma(i) \neq i$
and let $\gamma(a) = i$.
We make a few simple observations.
\begin{itemize}
\item We have $a > 1$.
\item We have $\gamma(a-1) > i$.
\item On the right of $i$, we have
$\gamma(a-1) - 1, \gamma(a-1) - 2, \cdots, i+1$.
\end{itemize} 
\end{rem}

\begin{thm}
\label{T: BPD to PD bijection}
The map $\flip$
is a weight-reversing bijection
between chains from $w$ to $w_0$ compatible with $(n-1, \cdots, 1)$
and chains from $w$ to $w_0$ compatible with $(1, \cdots, n-1)$.
Consequently, together with Theorem~\ref{T: BPD to chain} 
and~\cite{LeSo},
it can be viewed as a bijection from $\BPD(w)$ to $\PD(w)$.
\end{thm}
\begin{proof}
Let $C$ be a chain from $w$ to $w_0$ 
that is compatible with $(n-1, \cdots, 1)$.
We first check $\flip(C)$ 
is indeed a chain from $w$ to $w_0$
compatible with $(1, \cdots, n-1)$. 
More explicitly, 
we verify each $C_\gamma$ in 
Definition~\ref{D: BPD PD Bijection} is compatible with $\gamma$ recursively. 
Assume $C_\gamma = (w_1, \cdots, w_n)$ 
is compatible with $\gamma$ 
for some $\gamma$ that is not $[1, \cdots, n-1]$.
Let $i$ be the smallest with $\gamma(i) \neq i$
and suppose $\gamma(a) = i$.
The first two observations in Remark~\ref{R: gamma} say that $a > 1$ and $\gamma(a-1) > i$.
By our inductive hypothesis, 
we know $(w_{a-1}, w_a, w_{a+1},\cdots, \gamma(n-1))$
is compatible with 
$(\gamma(a-1), i,w_a, w_{a+1},\cdots, w_n)$.
By the third observation in Remark~\ref{R: gamma}, 
$i$ is the smallest number among
$i, \gamma(a+1),\cdots, \gamma(n-1)$
and all these numbers are smaller than $\gamma(a-1)$.
By Lemma~\ref{L: Stay in S^k (perm)},
$w_{a+1} \in S_n^{\gamma(a-1) \searrow}$.
By Lemma~\ref{L: stay in se k},
$w_{a-1} \in S_n^{\searrow i}$.
Now by Corollary~\ref{C: growth 2},
we know $\growth_{\gamma(a-1), i}(w_{a-1}, w_a, w_{a+1})$
will be compatible with $(i, \gamma(a-1))$.

It remains to sketch the inverse of $\flip$.
The description is similar to the 
description of $\flip$.
Let $C$ be a chain from $w$ to $w_0$
compatible with $(1, \cdots, n-1)$.
We will recursively compute chains $D_\gamma$,
where $\gamma$ range over certain permutations in $S_{n-1}$.
For the base case, we set $D_{[1, \cdots, n-1]} = C$.
Suppose we have $D_\gamma = (w_1, \cdots, w_n)$
where $\gamma$ is not $[n-1, \cdots, 1]$. 
Let $j$ be the largest such that 
$\gamma(j-1) < \gamma(j)$.
Let $\growth_{(\gamma(j-1), \gamma(j))}(w_{j-1}, w_j, w_{j+1}) 
= (w_{j-1}, w_j', w_{j+1})$.
Then define
$$D_{\gamma t_{j-1, j}} = (w_1, \cdots,w_{j-1}, w_j', w_{j+1},\cdots, w_n).$$
Eventually, we will obtain $D_{[n-1, \cdots, 1]}$.
We can similarly prove 
$D_{[n-1, \cdots, 1]}$ is a chain from $w$
to $w_0$ compatible with $[n-1, \cdots, 1]$.
By Corollary~\ref{C: growth 2},
it is the inverse of $\flip$.
\end{proof}

\begin{conj}
We conjecture this bijection agrees with the existing 
bijection of Gao and Huang~\cite{GH}.   
\end{conj}
This conjecture has been verified up to $S_7$.

\section{Hybrid Pipedreams and Bruhat chains}
\label{S: HPD}
As shown in Corollary~\ref{C: Bruhat formula},
for each choice of $\gamma \in S_{n-1}$,
there is a formula for $\fS_w$
involving chains in the Bruhat order.
The PD formula is identified with the chain
formula when $\gamma = [1,2, \cdots, n-1]$
and the BPD formula is identified with
the chain formula when $\gamma = [n-1, \cdots, 2, 1]$.
It would be natural to ask whether 
the other choices of $\gamma$ correspond 
to some combinatorial objects that compute
$\fS_w$.
We sketch how recent work of Knutson and Udell
provides an solution for $2^{n-2}$ such $\gamma$.

Knutson and Udell~\cite{KU} introduced 
\definition{hybrid pipe dreams} (HPD).
These are certain fillings of a $n \times n$ grid
with two pipes crossing at most once.
Each row can be either a 
``PD row'' or a ``BPD row''.
A PD row has a pipe entering from the left edge but not
right edge.
In a PD row,
the following tiles are are allowed
$$
\ptile, \quad \htile, \quad \jtile, \quad \rtile, \quad \bumptile, \quad \btile,
$$
while $\ptile$ and $\htile$ 
are called \definition{weighty tiles}.
A BPD row has a pipe entering from the right edge but not
the left edge.
In a BPD row,
the following tiles are are allowed
$$
\ptile, \quad \htile, \quad \qtile, \quad \Ltile, \quad \vtile, \quad \btile,
$$
while $\btile$ is weighty.
Let $\tau$ be a binary string
of length $n$ consisting of letters $P$ and $B$.
A HPD is called a \definition{$\tau$HPD} if 
the $i\textsuperscript{th}$ entry of $\tau$ indicates whether
row $i$ is a PD row or a BPD row. 
To associate a $\tau$HPD with a permutation,
we label the rows from $1, \cdots, n$ as follows:
Label the PD rows
by $1, \cdots, i$ from top to bottom,
and then label the BPD rows
by $i+1, \cdots, n$ from bottom to top.
Then we associate the $\tau$HPD with the permutation
$w$ where the pipe from a row
with label $i$
ends in column $w(i)$.
Let $\tau\HPD(w)$
be the set of all $\tau\HPD(w)$ associated
with $w$.
We let the weight of a HPD $H$,
denoted as $\wt(H)$,
be a sequence of $n-1$ numbers
where the $i\textsuperscript{th}$
entry is the number of weighty tiles
in the row with label $i$.
Clearly, HPDs recover the PDs when $\tau$
has only $P$s and 
recover the BPDs (but upside down) when $\tau$ 
has only $B$s.

\begin{exa}
\label{E: HPD}
Take $\tau = PBBPB$.
Let $H$ be the following $\tau\HPD$:

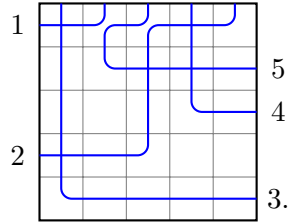
\[
\begin{tikzpicture}[x=1.5em,y=1.5em,thick,rounded corners, color = blue]
\draw[step=1,gray,thin] (0,0) grid (5,5);
\draw[color=black, thick, sharp corners] (0,0) rectangle (5,5);
\draw(5, 0.5)--(0.5,0.5)--(0.5,5);
\draw(0, 1.5)--(2.5,1.5)--(2.5,4.5)--(4.5,4.5)--(4.5,5);
\draw(5, 2.5)--(3.5,2.5)--(3.5,5);
\draw(5, 3.5)--(1.5,3.5)--(1.5,4.5)--(2.5,4.5)--(2.5,5);
\draw(0, 4.5)--(1.5,4.5)--(1.5,5);
\node[color=black] at (-0.5,4.5) {$1$};
\node[color=black] at (5.5,3.5) {$5$};
\node[color=black] at (5.5,2.5) {$4$};
\node[color=black] at (-0.5,1.5) {$2$};
\node[color=black] at (5.5,0.5) {$3$.};
\end{tikzpicture}
\]
We write the label of a PD (resp. BPD) row
on the left (resp. right) of the row. 
To find the permutation $w$ associated with $H$,
we trace the pipes. 
The pipe from the row with label $1$ goes
to column $2$, so we know $w(1) = 2$.
Following this procedure, 
we have $w = [2,5,1,4,3]$.
The heavy tiles are: 
$(1,1), (1,4), (3,2), (4,1),$ and $(4,2)$.
Thus, $\wt(H) = (2,2,0,1)$.
\end{exa}

\begin{rem}
\label{R: HPD double crossing}
Similar to Remark~\ref{R: BPD double crossing}, 
the condition ``two pipes
crossing at most once''
is equivalent to the
following requirement:
For each $\ptile$,
if the horizontal pipe
comes from a row labeled by $h$ 
and the vertical pipe comes 
from a row labeled by $v$,
we have $h > v$ if this $\ptile$
is in a BPD row and $h < v$ otherwise.   
\end{rem}

Knutson and Udell~\cite{KU} showed that 
the HPDs also compute the Schubert polynomials:
$$
\fS_{w} = \sum_{H \in \tau\HPD(w)} x^{\wt(w)},
$$
for any choice of $\tau$ and $w \in S_n$.
Since there are $2^n$ choices of $\tau$,
they obtain $2^n$ formulas of Schubert polynomials. 
We are going to show their formulas correspond
to $2^{n-2}$ Bruhat chain formulas.

The first step is to associate each HPD
with a Bruhat chain. 
Similar to the notion of $i$-row BPD, 
we introduce the following definition. 
For $i = 0, \cdots, n$,
we let $H_i$ be the filling obtained
by keeping the bottom $i$ rows of $H$
together with the row labels. 
We then define \definition{$i$-row HPDs}
as the set of all possible $H_i$
where $H$ ranges over all possible HPDs. 
We first derive two simple properties
of this notion.
\begin{lemma}
\label{L: HPD}
Let $H$ be a HPD. 
Assume the smallest row label in $H_i$
is $k$.
If we look at
the top edge of $H_i$,
there are no pipes
exiting from the $k-1$ columns on the right.
\end{lemma}
For instance, consider the $H$ in 
Example~\ref{E: HPD}.
The smallest label in $H_3$ is $2$.
Indeed, there are no pipes exiting from 
the rightmost column in $H_3$.
\begin{proof}
Prove by induction on $i$.
When $i = n$, 
we must have $k = 1$.
Our claim is trivial. 

Now suppose $i > 1$.
Assume toward contradiction that 
there is a pipe exiting from 
the top of $H_i$ in column $c > n - k +1$.
We compare $H_i$ and $H_{i-1}$.
Suppose the highest row of $H_{i-1}$
is a BPD row,
then the smallest row label
in $H_{i-1}$ is also $k$.
By our inductive hypothesis,
there are no pipe exiting on the right of column 
$n-k+1$ in $H_{i-1}$.
Thus, the topmost tile in 
column $c$ of $H_{i-1}$ is $\qtile$,
which means the pipe from the topmost row
exits on the right of column $c$.
Contradiction.
Now suppose the highest row of $H_{i-1}$
is a PD row, so it has label $k - 1$.
By our inductive hypothesis,
there are no pipe exiting 
on the right of column 
$n-k+2$ in $H_{i-1}$.
The topmost tile in 
column $c$ of $H_{i-1}$ is
$\ptile$, $\rtile$ or $\bumptile$.
In all cases, we know there
is a pipe exiting 
on the right of column 
$c$.
We have a contradiction since
$c \geq n - k + 1$.
\end{proof}

Next, we attach a permutation $u$ to an $i$-row HPD $H$ as follows.
For $k_1 \leq j \leq k_2$,
we let $u(j)$ be the number such that
the pipe from the label with label $i$
goes to column $u(j)$ in $H$.
By Lemma~\ref{L: HPD},
these numbers do not involve $n-k_1+2, \cdots, n$.
Then $u$ is determined if we require
$u \in S_n^{\searrow k_1} \cap S_n^{k_2 \searrow}$
(See Definition~\ref{D: se k} and 
Definition~\ref{D: k se}).

\begin{exa}
\label{E: i-row HPD}
Let $H$ be the HPD in Example~\ref{E: HPD}.
Then $H_5, \cdots, H_1,  H_0$ are associated 
with permutations 
$$
[2,5,1,4,3], \: [5,3,1,4,2], \: [5,3,1,4,2], 
\: [5,3,1,4,2], \: [5,4,1,3,2], \: [5,4,3,2,1].
$$
\end{exa}

We need one last lemma before we can
turn HPDs into Bruhat chains. 

\begin{lemma}
\label{L: HPD2}
Let $H$ be an $i$-row HPD
with topmost row labeled by $n$.
Let $H'$ be the $(i-1)$-row HPD
obtained by removing the topmost row of $H$.
Then $H$ and $H'$ are associated with the same permutation.
Moreover, $H$ is the only $i$-row
HPD whose highest row is labeled by $n$ 
and bottom $i-1$ rows agree with $H'$.
\end{lemma}
\begin{proof}
If $i = 1$, our lemma is trivial.
Otherwise, 
the highest row of $H$ must be a BPD row. 
Let $k$ be the smallest label in $H'$,
so labels in $H$ are $k, \cdots, n$.
There are $n - k + 1$ pipes 
exiting from the top of $H$.
By Lemma~\ref{L: HPD},
they are exiting from the first $n - k + 1$
columns. 
Thus, the tiling in this topmost 
row is uniquely determined by $H'$:
Let $c \leq n-k+1$ be the only number 
such that there is no pipe
exiting from the top of column $c$ in $H'$.
Then in the highest row of $H$,
the pipe from this row exits from column $c$.
and all other pipes go straight up.
This is equivalent to our lemma. 
\end{proof}

Now let $H$ be a $\tau$-HPD.
We may read the labels of rows
from top to bottom, 
obtaining a sequence of $n$ numbers. 
Let $\gamma_\tau$ be this sequence
without $n$.
Notice that $\gamma_\tau$
depends only on $\tau$.
Let $u_i$ be the permutation of 
$H_i$.
Assume row $k_0$ has label $n$.
By Lemma~\ref{L: HPD2},
$u_{k_0} = u_{k_0-1}$.
Define $\chain_\tau(H)$ as 
a chain of $n$ permutations 
obtained by removing $u_{k_0}$ from 
$(u_n, \cdots, u_0)$.
It is clear that $\chain_\tau$
recovers $\chain_\PD$ 
and $\chain_\BPD$ when $\tau$
has only $P$s or only
$B$s respectively.

\begin{exa}
Let $\tau = PBBPP$
and $H$ be the $\tau$HPD in Example~\ref{E: HPD}.
Then $\gamma_\tau = (1,4,2,3)$.
Following Example~\ref{E: i-row HPD},
we have
$$
\chain_\tau(H) = 
([2,5,1,4,3], \: [5,3,1,4,2], 
\: [5,3,1,4,2], \: [5,4,1,3,2], \: [5,4,3,2,1]).
$$
Notice that this chain is compatible with 
$(1,4,2,3)$.
Moreover, the length of the increasing $i$ 
chain agrees with
$n - i$ minus the number of weighty tiles
in the row with label $i$ in $H$.
\end{exa} 

\begin{thm}
\label{T: HPD to chain}
The map $\chain_\tau(\cdot)$
is a bijection from
$\tau$-HPDs of $w$
to chains from $w$ to $w_0$
compatible with $\gamma_\tau$.
Moreover, if the row with label $k$ of 
a $\tau$-HPD $H$ has $m$ weighty tiles, 
then the increasing $k$-chain
in $\chain_\tau(H)$
has length $n - k - m$.
\end{thm}

Our proof will mimic the proof of Theorem~\ref{T: BPD to chain}. 
We first check $\chain_\tau(H)$ 
is indeed compatible with $\gamma_\tau$,
which equivalent to checking the following:
\begin{lemma}
\label{L: HPD one row to chain}
Let $H$ be an $i$-row HPD with $i \geq 1$.
Let $H'$ be the $(i-1)$-row HPD obtained 
by removing the highest row of $H$.
Assume $H$ is associated with $u$
and $H'$ is associated with $w$.
Let $k$ be the row label of 
the highest row in $H$.
If $k \neq n$,
we know there is an increasing $k$-chain
from $u$ to $w$.
Moreover, a number $j$ is not among $w(1), \cdots, w(k)$ and is fixed in the $k$-chain
if and only if $(1, j)$ in $H$ is weighty.
\end{lemma}
\begin{proof}
First, suppose row $1$ of $H$ is a BPD row. 
Let $k_1$ be the smallest label in $H$.
We know $u, w \in S_n^{\searrow k_1}$.
After ignoring the first $k_1 - 1$ numbers 
in both $u$ and $w$, 
we obtain permutations in $S_{n - k_1 + 1}$.
In addition, by Lemma~\ref{L: HPD},
the rightmost $k_1 - 1$ tiles 
in row $1$ of $H$ must be 
$\htile$.
Then the lemma can be proved
similarly as Lemma~\ref{l: One row to chain}.

Now we assume row $1$ of $H$ is a PD row. 
Let $k_2$ be the largest row label
in $H'$.
We know $w \in S_n^{k_1 + 1 \searrow} \cap S_n^{k_2 \searrow}$.
Let $j_1 < \cdots < j_d$ be the column indices
of non-heavy tiles in the first $n - k$ columns of row $1$ of $H$.
Then we claim the increasing $k$-chain
from $u$ to $w$ can be expressed as:
$$
w_1 \xrightarrow{t_{k, w^{-1}(j_1)}}
\cdots \xrightarrow{t_{k, w^{-1}(j_d)}} w_{d+1} = w.
$$
In other words, 
$w_i$ is obtained from $w_{i+1}$ 
by swapping the value 
$j_i$ with the 
$k\textsuperscript{th}$ entry. 
We first check this is indeed a 
$k$-chain.
Take $i \in [d]$.
We know $j_i \leq n-k$.
Since $w \in S_n^{\searrow k+1}$,
the first $k$ numbers in $w$
are $n, \cdots, n- k + 1$.
Thus, $j_i$ is not among the first
$k$ numbers of $w$, or $w_{i+1}$.
Let $j_{d+1} = n - k + 1$ by convention.
Then we know $w_{i+1}(k) = j_{i+1}$.
For any $j_d < c < j_{d+1}$,
we know $(1, c)$ in $H$ is heavy,
so it is either $\ptile$ or $\htile$.
The horizontal pipe in it
enters the row from the bottom 
of $(1, j_i)$,
so it is from a row with label
$w^{-1}(j_i)$.
If $(1, c)$ is $\ptile$ with the vertical pipe coming from 
a row with label $v$,
we know $v > w^{-1}(j_i)$ 
by Remark~\ref{R: HPD double crossing}.
In other words, 
$c$ appears on the right of $j_i$
in $w$ and stays unswapped in
the chain, so $c$ is on the 
right of $j_i$ in $w_{i+1}$.
If $(1, c)$ is $\htile$,
we know $c$ is not among 
the first $k_2$ numbers in $w$
but $j_i$ is.
Again, we conclude $c$ is on the 
right of $j_i$ in $w_{i+1}$.
In either case, 
$w_i \lessdot_k w_{i+1}$.
Then the chain is an increasing
$k$-chain since the smaller
numbers being swapped
are $j_1 < \cdots < j_d$.
Numbers not among the first $k$
numbers in $w$ are $1, \cdots, n - k$.
Such a number $c$ is unswapped in the 
chain if and only if $(1, c)$ in $H$ 
is weighty.

It remains to check $w_1 = u$.
We have $w_1 \in S_n^{\searrow k}$.
By Lemma~\ref{L: Stay in S^k (perm)},
$w_1 \in S_n^{k_2 \searrow}$.
We just need to check 
$w_1(i) = u(i)$ for each $k \leq i \leq k_2$.
Notice that $w_1$ is obtained from
$w$ by applying the following 
changes: 
Replace $j_i$ by $j_{i+1}$ 
for $i \in [d]$
and then let set the 
$k\textsuperscript{th}$ number to be $j_1$.
We check $u$ is obtained from $w$
in the same way. 
\begin{itemize}
\item The first $j_1 - 1$ tiles
in row $1$ of $H$ are 
$\ptile$ or $\htile$.
Thus, the pipe from row the row 
with label $k$ 
goes to column $j_1$. 
\item Suppose a pipe enters on
the bottom of $(1, j_i)$.
Since this tile is not weighty,
this pipe exits the tile
from the right.
There are only weighty tiles between
$(1, j_i)$ and $(1, j_{i+1})$,
so this pipe enters $(1, j_{i+1})$
from the left. 
Then since $(1, j_{i+1})$ 
is not weighty,
it exits from the top of this tile. 
Thus, the $j_{i}$ in $w$ 
is replaced by $j_{i+1}$ in $u$.
\item Suppose a pipe enters on
the bottom of $(1, c)$ which is weighty.
Then this tile is $\ptile$
and it exits from the top.
Correspondingly, 
the value $c$ in $w$ is not
changed in $u$.
\end{itemize}
Thus, we have $w_1 = u$.
\end{proof}

This Lemma allows us
to see $\chain_\tau(H)$
is compatible with $\gamma_{\tau}$
when $H$ is a $\tau$HPD.
It remains to construct the 
inverse of $\chain_\tau$.
To do so, 
we need to characterize what 
permutations can be associated 
with a $k$-row HPD.

\begin{lemma}
\label{L: Stay in S^k (perm)2}
Take $w \in S_n^{\searrow k + 1}$.
If there is an increasing $k$-chain
from $u$ to $w$,
then $u \in S_n^{\searrow k}$.
In other words, 
the chain keeps swapping the 
$k\textsuperscript{th}$ entry
with some number on its right.
\end{lemma}
\begin{proof}
Let the increasing 
$k$-chain from $u$ to $w$ be
$(w_1, \cdots, w_d)$.
Consider how to obtain $w_{d-1}$ from $w_d$.
Since $w_d = w \in S_n^{\searrow k + 1}$,
we know $w_{d-1}$ is obtained from $w_d$
by swapping the $k\textsuperscript{th}$
entry with someone on its right.
Then since the chain is increasing,
the first $k-1$ entries stays fixed 
throughout the chain,
so $u \in S_n^{\searrow k}$.
\end{proof}

Then constructing the inverse of 
$\chain_\tau$ relies on the following:
\begin{lemma}
\label{L: HPD add a row}
Let $H$ be an $(i-1)$-row HPD
with row labels $k_1, \cdots, k_2$
and is associated with $w$.
\begin{itemize}
\item Suppose there is an increasing
$k_2 + 1$ chain from $u$ to $w$.
Then there exists a $i$-row HPD
associated with $u$ 
such that its bottom-most $i-1$
rows agree with $H$
and top-most row is a BPD row. 
\item Suppose there is an increasing
$k_1 - 1$ chain from $u$ to $w$.
Then there exists a $i$-row HPD
associated with $u$ 
such that its bottom-most $i-1$
rows agree with $H$
and top-most row is a PD row.     
\end{itemize}
\end{lemma}
\begin{proof}
For the first statement,
since $w \in S_n^{\searrow k_1}$,
we know $u \in S_n^{\searrow k_1}$.
After ignoring the first $k_1 - 1$
numbers in $w$ and $u$,
we may view them as permutations
from $S_{n - k_1 + 1}$.
We may construct a row of tilings
with $n - k_1 + 1$ tiles
as in Lemma~\ref{L: one row BPD2}
(but upside down).
Then we simply append $(k_1 - 1)$
$\htile$s on its left to obtain
a row of $n$ tiles.
We place them on top of $H$
to obtain the desired $i$-row HPD.

Now we consider the second statement.
We know $w \in S_n^{\searrow k_1 + 1}$.
By Lemma~\ref{L: Stay in S^k (perm)2},
$u$ is obtained from $w$
by swapping the values $j_d > \cdots > j_1$ with the $k_1\textsuperscript{th}$
entry. 
If we let $j_{d+1} = n - k_1 + 1$,
then $u$ is obtained from $w$
by swapping $j_{i+1}$
with $j_i$ for $i = d,\cdots, 1$.

Now we construct a row of $n$ tiles 
using tiles of a PD row.
The rightmost $k-1$ tiles are $\btile$.
Take $c \in [n - k + 1]$ 
that is not one of $j_1, \cdots, j_{d+1}$.
We let $(1, c)$ be $\htile$
if $w^{-1}(c) > k_2$ and $\ptile$ otherwise. 
Then $(1,j_1), \cdots, (1, j_{d+1})$ can be 
$\btile$, $\rtile$, $\jtile$ 
or $\bumptile$.
We may determine $(1, j_m)$
for $m = 1,\cdots, d+1$ using the following table:
\[
\begin{array}{c|c|c}
& \text{If $m = 1$ or $(1, j_m - 1)$
has a pipe exiting from the right}& \text{otherwise} \\
\hline
\text{$k_1 < w^{-1}(j_m) \leq k_2$} 
& \bumptile & \rtile\\
\hline
\text{Otherwise}
& \jtile & \btile\\
\end{array}
\]

From our construction, it is clear that 
$(1, c)$ is weighty if and only if
$c < n-k_1 + 1$ and the value 
$c$ is fixed throughout the increasing $k_1$-chain.

First, we need to show the pipes are connected.
In other words, we need to check
for each tile $(1,c)$, it has a pipe coming from 
the left edge if and only if there is a pipe
exiting from the right edge of $(1, c-1)$ or $c = 1$.
We prove this claim by considering three cases of $c$:
\begin{itemize}
\item If $c > n - k_1 + 1$, 
we know $(1, c)$ is $\btile$.
We just need to make sure there are
no pipe exiting from the left of $(1, n-k_1+1)$.
We know $j_{d+1} = n-k+1$
and $w(k_1) = n - k_1 + 1$.
By our construction,
$(1, n - k_1 + 1)$ is $\btile$
or $\jtile$.
\item If $c = j_i$ for some $i$,
this statement follows from our construction.
\item Otherwise, $c \leq n - k + 1$
and $(1, c)$ is weighty, 
so it is $\ptile$ or $\htile$.
We need to verify $(1, c-1)$ has a pipe
exiting from the right edge if $c > 1$.
This is immediate if $(1, c)$ is also weighty.
If not, 
we know the increasing $k_1$-chain
fixes the value $c$ but swaps the value 
$c - 1$ with the $k_1\textsuperscript{th}$
entry. 
Thus, $c-1$ is on the left of $c$ in $w$.
Since $w \in S_n^{k_2 \searrow}$,
we know $c-1$ is among the first $k_2$ numbers
in $w$.
Thus, $(1, c-1)$ is $\bumptile$ or $\rtile$
by construction.  
\end{itemize}

Now we append this row on top of $H$.
We need to show what we get is a $i$-row HPD.
In other words, we need to check the
bottom edge of the row we constructed
is consistent with the top edge of $H$.
The $(1,c)$ of $H$ 
has a pipe exiting from the top
if and only if $k_1 < w^{-1}(c) \leq k_2$.
By our construction, 
this happens if and only if
the $(1,c)$ we construct
is $\ptile$, $\rtile$ or $\bumptile$,
which is equivalent to:
the $(1,c)$ we construct
has a pipe coming from the bottom edge.

Finally, assume this $i$-row HPD we obtained 
is associated with the permutation $u'$.
We need to verify $u = u'$.
We know 
$u, u' \in S_n^{\searrow k_1} \cap S_n^{k_2 \searrow}$,
so we need to make sure $u(t) = u'(t)$
for $k_1 \leq t \leq k_2$.
Recall that $u$ is obtained from $w$
by changing $j_m$ into $j_{m+1}$
for $m \in [d]$ and setting
the $k_1\textsuperscript{th}$
number into $j_1$.
We verify $u'$ can be obtained from
$w$ in the same way.
\begin{itemize}
\item The row $1$ in our $i$-row HPD
has label $k_1$.
A pipe enters from the left edge of this row.
By our construction $(1,j_1)$ is $\jtile$ or $\bumptile$
and all tiles on its left are $\ptile$ or $\htile$.
Thus, the pipe from row $1$ exits 
from column $j_1$,
so $u'(k_1) = j_1$.
\item Take $k_1 < t \leq k_2$
such that $w(t) = j_m$ for some $m \in [d]$.
Then the tile $(1,j_m)$ 
in our constructed row
is $\rtile$ or $\bumptile$.
Tiles between $(1, j_m)$
and $(1, j_{m+1})$
are $\ptile$ or $\htile$.
Then the tile $(1, j_{m+1})$
is $\jtile$ or $\bumptile$.
Thus, the pipe entering 
from the bottom of $(1,j_m)$
exits from the top of $(1, j_{m+1})$.
In other words,
the $j_m$ in $w$ is replaced 
by $j_{m+1}$ in $u'$.
\item Take $k_1 < t \leq k_2$
such that $w(t)$ is not among
$j_1, \cdots, j_{d+1}$.
Then by construction,
$(1, w(t))$ of our constructed
row is $\ptile$. 
The pipe entering 
from the bottom of $(1,w(t))$
exits from the top of $(1, w(t))$.
In other words, 
the value $w(t)$ stays 
the unchanged in $u'$.
\qedhere
\end{itemize}
\end{proof}

Finally, we prove Theorem~\ref{T: HPD to chain}.
We omit the details since the proof
is essentially the same as
the proof of Theorem~\ref{T: BPD to chain}.

\begin{proof}[Proof of Theorem~\ref{T: HPD to chain}]
Take $H \in \tau\HPD(w)$.
By Lemma~\ref{L: HPD one row to chain}, we know $\chain_\tau(H)$ is a chain
from $w$ to $w_0$ compatible with
$\gamma_\tau$
and respects the weight. 

Then we show $\chain_\tau$ is injective.
Suppose $\chain_\tau(H') = \chain_\BPD(H)$
where $H'$ is another $\tau$HPD.
Then for each $k \in [n-1]$,
the starting points and destinations of each pipe
are the same in the row with label $k$
in $H$ and $H'$.
The row with label $n$ agree
by Lemma~\ref{L: HPD2}.
Finally, a right inverse of $\chain_\tau$
can be constructed by repetitively
invoking Lemma~\ref{L: HPD add a row}.
Since $\chain_\BPD(\cdot)$ is injective,
this map is its inverse. 
\end{proof}

\begin{rem}
Clearly there are $2^{n-2}$ distinct $\gamma_\tau$
when $\tau$ ranges over all $2^n$ possible choices. 
Therefore, 
the HPDs correspond to $2^{n-2}$ formulas
of those $(n-1)!$ chain formulas from Corollary~\ref{C: Bruhat formula}.
\end{rem}

\section{Future Directions}
\label{S: Future}
In this paper,
we realize that
PDs and BPDs correspond to two ``extremes'' in the family
of the Bruhat chains formulas:
There are $(n-1)!$ different choices of the $\gamma$
in Corollary~\ref{C: Bruhat formula}
where PDs and BPDs correspond to 
the choice $\gamma = (1, \cdots, n-1)$ and 
$\gamma = (n-1, \cdots, 1)$ respectively.
It would be natural to ask if the rich mathematics
of PDs and BPDs 
extend to the whole family. 
\begin{itemize}
\item Monk's rule~\eqref{eq: Monk} can been proved 
by building bijections between combinatorial objects
that compute the Schubert polynomials.
Bergeron and Billey~\cite{BB} gave such a proof on PDs
and Huang~\cite{Huang_Monk} gave such a proof on BPDs.
Their proofs involve very different combinatorial
moves 
but sharing some similar intrinsic ideas,
such as allowing two pipes to cross twice temporarily. 
Can one give a bijective proof of the Monk's rule for 
arbitrary $\gamma$?
\item Kashiwara~\cite{Kas90, Kas91} introduced crystal 
operators. 
Following~\cite{BS_crystal},
they may be defined as operators on combinatorial objects,
such as words or tableaux, 
satisfying certain axioms.
There is a recent surge on defining crystal operators on PDs~\cite{AS,GMS,MS_crystal}.
One application of defining crystal operators on
PDs is to interpret the Schubert-to-key
expansion given by Reiner and Shimozono~\cite{RS}:
PDs that are connected via the crystal moves 
in $\PD(w)$ form a key polynomial
(see~\cite{GMS} for more details).

In work under preparation, 
Huang, Shimozono and Yu define crystal moves
on BPDs. 
They will show that the Gao-Huang bijection
gives a crystal isomorphism.
The crystal operators on PDs and BPDs
are achieved via different local
moves: ``chute moves'' on PDs and ``droop moves'' on BPDs. 
Is there a unifying way to define crystal operators
on the chains in Corollary~\ref{C: Bruhat formula}
for arbitrary choice of $\gamma$?

\item 
Edelman and Greene~\cite{EG} defined 
an insertion algorithm that sends $P \in \PD(w)$ to a pair of tableaux: 
an insertion tableau and a recording tableau. 
The shape of the insertion tableaux index the 
Schur expansion of the Stanley symmetric 
function of $w$~\cite{Stan}.
The recording tableaux are semi-standard Young tableaux (SSYTs).

In work under preparation, 
Huang, Shimozono and Yu refine
an insertion algorithm of~\cite{LLS} to give a bijection
that sends BPDs to a pair 
where the first entry is a ``EG BPD'' and the 
second entry is a SSYT. 
The EG BPDs are certain BPDs defined in~\cite{LLS}, also indexing the Schur expansion 
of Stanley symmetric functions.
Is there a unifying way
to define an insertion algorithm
that sends chains in Corollary~\ref{C: Bruhat formula}
for arbitrary $\gamma$
to a pair of combinatorial objects,
where the first object index the Stanley-to-Schur
expansion and the second object is a SSYT?

\item Grothendieck polynomials, introduced by Lascoux and Sch\"utzenberger~\cite{LS:Groth},
are K-theoretic analogues of Schubert polynomials. 
They are inhomogeneous polynomials whose
bottom-degree homogeneous components are Schubert polynomials. 
They can be computed using extensions of PDs~\cite{FK_Gro} 
and BPDs~\cite{W} where two pipes are allowed 
to cross more than once.
Is there a way to extend the notion of chains
in Corollary~\ref{C: Bruhat formula}, 
obtaining $(n-1)!$ different chain formulas
for Grothendieck polynomials
that naturally recover both PD and BPD formulas. 

\item 
The \definition{double Schubert polynomial} $\fS_w(\textbf{x},\textbf{y})$
is in $x_1, \cdots, x_{n-1}$
and $y_1, \cdots, y_{n-1}$.
It recovers $\fS_w$ after setting
each $y_i$ to $0$
and can be computed using 
PDs and BPDs:
For $P \in \PD(w)$ (resp. $\BPD(w)$),
let $\WT(P)$ 
be the product over $\ptile$ (resp. $\btile$) in $P$,
where the tile in row $i$ column $j$ gives $(x_i - y_j)$.
By~\cite{KM, W}, 
$$\fS_w (\textbf{x},\textbf{y}) = \sum_{P \in \PD(w)} \WT(P)
= \sum_{P \in \BPD(w)} \WT(P).$$

Recently, Knutson and Udell showed  
HPDs can also compute double Schubert polynomials. 
For an arbitrary $\tau$, 
let $H$ be a $\tau$-HPD associated with $w$.
Let $\WT(H)$ 
be the product over weighty tiles in $H$,
where the tile in the row with label $i$ 
and column $j$ gives $(x_i - y_j)$.
By~\cite{KU}, 
$\fS_w (\textbf{x},\textbf{y}) = \sum_{H \in \tau-\HPD(w)} \WT(H)$.
In other words, 
for $2^{n-2}$ choices of $\gamma$,
one can compute $\fS_w (\textbf{x},\textbf{y})$
using the chains compatible with $\gamma$.
How can
the $(n-1)!$ different chain formulas in Corollary~\ref{C: Bruhat formula} 
be augmented to compute 
$\fS_w (\textbf{x},\textbf{y})$?
\end{itemize}

We would like to give a conjecture for the last question above.
Take $\gamma \in S_{n-1}$
and let $C = (w_1, \cdots, w_n)$ 
be a chain compatible with $\gamma$.
Define $\WT_\gamma(C)$
as $\prod_{i = 1}^{n-1} \prod_t (x_{\gamma_i} - y_{w_i(t)})$,
where $t$ runs over all 
$t > \gamma_i$
such that $w_i(t) = w_{i+1}(t)$.
After setting all $y_i$ to $0$, $\WT_\gamma(C)$
recovers $x^{(n-1, \cdots, 1) - \gamma^{-1}(\wt(C))}$.
The following conjecture
extends Corollary~\ref{C: Bruhat formula}
and has been checked for all $w \in S_n$ for $n \leq 8$ and all $\gamma \in S_{n-1}$:
\begin{conj}
For $\gamma \in S_{n-1}$, 
we have $\fS_w$ $(\textbf{x},\textbf{y}) = \sum_{C: \textrm{chain from $w$ to $w_0$ compatible with $\gamma$}} \WT_\gamma(C)$. 
\end{conj}

This conjecture agrees with the PD (resp. BPD) formula
when $\gamma = [1, \cdots, n-1]$ 
(resp. $[n-1, \cdots, 1]$)
via the bijections $\chain_\PD$ (resp. $\chain_\BPD$) by Remark~\ref{R: Double PD} (resp. Remark~\ref{R: Double BPD}).

Finally, we outline two possible strategies
to approach the questions above and
identify one major challenge in each strategy. 
One may first focus on the choices of $\gamma$ 
that correspond to HPDs. 
The pipe formulas allow one to think about local moves
more intuitively. 
After this step, 
it would be challenging to extend the arguments to general
$\gamma$ without ``pipe-like'' formulas. 
Therefore, to take this approach, 
it is crucial to come up with combinatorial formulas
involving pipes that naturally
model the chain formulas for each possible $\gamma$.

One may also start from known results on PDs or BPDs 
and gradually move to general $\gamma$ by swapping 
two adjacent numbers in $\gamma$ each time.
For instance, by Corollary~\ref{C: growth 1},
we have a bijection between chains from $w$ to $w_0$
that are compatible with $[1, 2, 3, \cdots, n-1]$
and $[2,1,3, \cdots, n-1]$ using Lenart's growth diagram.
Crystal operators are already defined on the former chains. 
To make this bijection a crystal isomorphism, 
there is a unique way to define crystal operators
on the latter chains. 
Then one can attempt to find a simpler description
of these operators.
This approach relies on bijections 
between chains of slightly different $\gamma$
and Corollary~\ref{C: growth 1} does not provide
enough such bijections.
Thus, it is crucial to resolve 
Problem~\ref{prob: chain} in more cases.
Unfortunately, Lenart's growth diagram does 
not do this in general, 
so one has to come up with new techniques
and operations on Bruhat chains. 

\section{Acknowledgment}
We thank Yibo Gao for suggesting 
the problem of finding a BPD analogue 
of the Fomin-Stanley construction.
We thank Yibo Gao and Zachary Hamaker 
for many important suggestions, 
including using Lenart's Growth diagram to
obtain a bijection between PDs and BPDs.
We thank Yibo Gao, Yuhan Jiang and Brendon Rhoades for
carefully reading an earlier version
of this paper and providing valuable comments.

\bibliographystyle{alpha}
\bibliography{main}{}
\end{document}